\documentclass[reqno]{amsart}
\usepackage{fullpage}
\usepackage{array}
\usepackage{mathrsfs}
\usepackage{multirow}
\usepackage{hyperref}
\usepackage{graphicx,color}
\usepackage{float}
\usepackage{ytableau}
\usepackage{amsmath}
\usepackage{amssymb}

\newtheorem{theorem}{Theorem}[section]
\newtheorem{example}[theorem]{Example}
\newtheorem{lemma}[theorem]{Lemma}
\newtheorem{corollary}[theorem]{Corollary}
\newtheorem{remark}[theorem]{Remark}
\newtheorem{proposition}[theorem]{Proposition}
\theoremstyle{definition}

\newcommand{\var}{\operatorname{var}}
\newcommand{\Id}[1]{\operatorname{Id}(#1)}

\usepackage{comment}

\usepackage{xcolor}

\newcommand{\red}{\color{red}}
\newcommand{\blue}{\color{blue}}
\newcommand{\black}{\color{black}}

\begin{document}

\title{{Varieties of associative algebras \\ with slow central codimension growth}}



\author{
Wesley Quaresma Cota$^{a}$, Antonio Ioppolo$^{b,*}$ and Ana Cristina Vieira$^c$ 
}

\thanks{{\it E-mail addresses:} wesleyqc@ime.usp.br (Cota), antonio.ioppolo@univaq.it (Ioppolo) and  anacris@ufmg.br (Vieira).}

\dedicatory{$^{a}$ IME, USP, Rua do Matão 1010, 05508-090, São Paulo, Brazil \\ $^b$ DISIM, Universit\`a degli Studi dell'Aquila, Via Vetoio 1, 67100, L'Aquila, Italy \\ 
$^c$ ICEx, UFMG, Avenida Antonio Carlos 6627,
31123-970, Belo Horizonte, Brazil}

\thanks{\footnotesize W. Q. C.~is supported by FAPESP, grant nos.~2024/14914-9, 2025/05699-0 and 2026/06588-0.}
\thanks{A. I. is supported by GNSAGA-INdAM.}
\thanks{A. C. V. is supported by CNPq.}
\thanks{$^{*}$ Corresponding author}



\begin{abstract}     
Central codimensions measure the dimensions of spaces of multilinear polynomials modulo those taking values in the center of an algebra. We classify associative algebras over a field of characteristic zero whose central codimensions grow at most linearly. The classification is given up to equality of polynomial identities, in terms of explicit generating algebras and an equivalent finite exclusion criterion. We also give a direct structural characterization of the bounded case. As a consequence, we determine the five varieties of quadratic central codimension growth whose proper subvarieties all have at most linear central growth. As an application of the bounded-growth classification, we characterize algebras graded by an arbitrary finite group whose graded central codimensions are bounded. Up to graded PI-equivalence, these are direct sums of commutative and nilpotent graded algebras, possibly together with the Grassmann algebra on two generators endowed with the trivial grading.

\end{abstract}

\maketitle

\section{Introduction}\label{sec:introduction}

The codimensions $c_n(A)$, $n\geq 1$, of an associative algebra $A$ over a field $F$ of characteristic zero record the dimensions of the spaces of multilinear polynomials in $n$ variables modulo the identities of $A$. Their asymptotic behavior relates the structure of associative algebras to representation theory and asymptotic analysis. Two fundamental results underlie this study: Regev's theorem that the codimensions of a PI-algebra are exponentially bounded~\cite{RG}, and Kemer's dichotomy between polynomially bounded and exponential growth~\cite{Kemer}.

For every nonnilpotent PI-algebra, Giambruno and Zaicev~\cite{GiambrunoZaicev} proved that the exponential growth rate $d=\exp(A)$ exists and is a positive integer. A finer description is predicted by Regev's asymptotic conjecture,
$$
c_n(A)\sim\alpha n^t d^n,
$$
with $\alpha>0$, which separates the exponential rate $d$ from the polynomial factor $n^t$. This conjecture was established for unital PI-algebras by Berele and Berele--Regev~\cite{Berele,BR}.

For such nonnilpotent algebras, without assuming a unit, Giambruno and Zaicev~\cite{GZ14} proved that the polynomial exponent
$$
t=\lim_{n\to\infty}\frac{\log\bigl(c_n(A)/d^n\bigr)}{\log n}
$$
exists and takes an integer or half-integer value. For nonnilpotent algebras of polynomial codimension growth, $d=1$ carries no information about the polynomial degree. The invariant $t$ distinguishes growth at this finer level and motivates the structural problem of classifying algebras according to their polynomial degree~\cite{GLa,CMQ,Sandra,KoshlukovLaMattina2015}.

We study the corresponding problem for the central codimensions $c_n^z(A)$, which measure multilinear polynomials modulo those taking values in the center of $A$. Giambruno and Zaicev~\cite{GiZai2} established the existence and integrality of the central exponent and its equality with $d$ when $d\geq 2$. Together with Kemer's dichotomy, this implies that ordinary codimensions are polynomially bounded if and only if central codimensions are. Their polynomial degrees may nevertheless differ: the Grassmann algebra $\mathcal{G}_2$ on two generators satisfies $c_n(\mathcal{G}_2)=1+\binom{n}{2}$, whereas $c_n^z(\mathcal{G}_2)=1$.

Central polynomials form an endomorphism-invariant vector space, generally not an ideal, and low central codimension growth is compatible with nonzero central products of radical components. The structural issue is to control the interactions among these components. We do so through a finite family of excluded algebras and reductions to explicit models up to PI-equivalence.

We classify, up to PI-equivalence, algebras with bounded or at most linear central codimensions, without assuming finite dimensionality or the existence of an identity element. Here PI-equivalence means equality of polynomial identities. In the bounded case, the only possibilities are $N$, $F\oplus N$ and $\mathcal{G}_2\oplus N$, where $N$ is nilpotent; hence the central codimensions are eventually identically zero or identically one. This classification also follows from~\cite[Theorem~5.6]{NascimentoCota}; here we give a direct structural proof and an explicit exclusion criterion. Our main result treats at most linear growth: every such algebra is PI-equivalent to a nilpotent algebra or to a finite direct sum of members of a fixed list of six finite-dimensional algebras, together with a nilpotent summand.

As a consequence, we determine the five $z$-minimal varieties of quadratic central codimension growth, that is, those whose proper subvarieties all have at most linear central codimension growth. In the ordinary theory, the generators of the minimal quadratic varieties~\cite{Sandra} form part of the building blocks used in the complete classification of varieties of quadratic codimension growth~\cite{CMQ}. This structural role motivates the study of the corresponding minimal varieties in the central setting.

As an application, we use the bounded-growth classification to characterize algebras with bounded graded central codimensions for gradings by any finite group. Gradings already play a structural role in Kemer's representability theorem~\cite{Kemer}: every PI-algebra has the same polynomial identities as the Grassmann envelope of a finite-dimensional $\mathbb{Z}_2$-graded algebra.

In the graded setting, a commutative algebra may have exponentially growing graded codimensions while all its graded central codimensions vanish. We prove that a graded algebra has bounded graded central codimensions if and only if it is graded PI-equivalent to $C\oplus N$ or $\mathcal{G}_2\oplus C\oplus N$, where $C$ is a commutative graded algebra, either zero or unital, $N$ is a nilpotent graded algebra, and $\mathcal{G}_2$ is trivially graded. The graded representability theorem of Aljadeff and Kanel-Belov~\cite{AljadeffKanelBelov2010} allows us to obtain this classification without assuming finite dimensionality.

\section{Polynomial identities and central polynomials}
\label{sec:preliminaries}

Throughout the paper, $F$ denotes a field of characteristic zero, and
all algebras are associative $F$-algebras, not necessarily with an
identity element. Let $F\langle X\rangle$ be the free nonunital associative $F$-algebra on a countable set $X=\{x_1,x_2,\ldots\}$ of variables. Thus, its elements are noncommutative polynomials with zero constant term.

A polynomial $f=f(x_1,\ldots,x_m)\in F\langle X\rangle$ is a
\emph{polynomial identity} of an algebra $A$ if
$f(a_1,\ldots,a_m)=0$ for every $a_1,\ldots,a_m\in A$. The set of
all polynomial identities of $A$ is denoted by
$\operatorname{Id}(A)$. It is a \emph{$T$-ideal}, that is, an ideal
of $F\langle X\rangle$ invariant under all endomorphisms of the
free algebra. The algebra $A$ is called a \emph{PI-algebra} if
$\operatorname{Id}(A)\neq\{0\}$. Given polynomials
$f_1,\ldots,f_r$, we write
$\langle f_1,\ldots,f_r\rangle_T$ for the $T$-ideal they generate.

Two algebras $A$ and $B$ are \emph{PI-equivalent}, written
$A\sim_T B$, if $\operatorname{Id}(A)=\operatorname{Id}(B)$.
The variety generated by $A$, denoted by $\operatorname{var}(A)$, is
the class of algebras $B$ such that $\textnormal{Id}(A)\subseteq \textnormal{Id}(B).$ In particular, subalgebras and homomorphic images of $A$ belong
to $\operatorname{var}(A)$.

We use the notation $[x_1,x_2]=x_1x_2-x_2x_1$ for the commutator, and longer commutators are left-normed:
\[
[x_1,\ldots,x_r]
=
[[x_1,\ldots,x_{r-1}],x_r],
\qquad r\geq 3.
\]

Let $$Z(A)=\{a\in A\mid ab=ba\text{ for all }b\in A\}$$ be the
center of $A$. A polynomial $f=f(x_1,\ldots,x_m)$ is called a
\emph{central polynomial} of $A$ if every evaluation of $f$ in
$A$ belongs to $Z(A)$. Following the terminology of
\cite{GiZai2}, identities are included among central polynomials;
a central polynomial that is not an identity is called a
\emph{proper central polynomial}. We denote the space of central
polynomials of $A$ by $\operatorname{Id}^z(A)$.

Unlike $\operatorname{Id}(A)$, the space $\operatorname{Id}^z(A)$
need not be an ideal. It is, however, a \emph{$T$-space}, namely,
an $F$-subspace of $F\langle X\rangle$ invariant under all
endomorphisms. When needed, we write
$\langle f_1,\ldots,f_r\rangle^T$ for the $T$-space generated
by $f_1,\ldots,f_r$. Clearly,
$\operatorname{Id}(A)\subseteq\operatorname{Id}^z(A).$ Moreover, for a variable $x_{m+1}$ not occurring in $f$, we have
$f\in\operatorname{Id}^z(A)$ if and only if $
[f,x_{m+1}]\in\operatorname{Id}(A).$ Consequently, the $T$-space of central polynomials is determined by the $T$-ideal of identities. In particular, if $A\sim_T B$ then 
$\operatorname{Id}^z(A)=\operatorname{Id}^z(B).$

Since the characteristic is zero, the multilinearization process allows us to
recover both $\operatorname{Id}(A)$ and $\operatorname{Id}^z(A)$ from their multilinear elements. For $n\geq 1$, let $S_n$ denote the symmetric group on
$\{1,\ldots,n\}$, and let
\[
P_n=\operatorname{span}_F
\{x_{\sigma(1)}\cdots x_{\sigma(n)}\mid \sigma\in S_n\}
\]
be the space of multilinear polynomials in $x_1,\ldots,x_n$.
Consider the quotient spaces
\[
P_n(A)
=
\frac{P_n}{P_n\cap\operatorname{Id}(A)}
\quad \mbox{ and }\quad 
P_n^z(A)
=
\frac{P_n}{P_n\cap\operatorname{Id}^z(A)}.
\]
Their dimensions
\[
c_n(A)=\dim_F P_n(A)
\quad \mbox{ and }\quad 
c_n^z(A)=\dim_F P_n^z(A)
\]
are the \emph{$n$th ordinary codimension} and the
\emph{$n$th central codimension} of $A$, respectively.

For comparison, the multilinear proper central polynomials are
measured by
\[
\Delta_n(A)
=
\frac{P_n\cap\operatorname{Id}^z(A)}
     {P_n\cap\operatorname{Id}(A)},
\qquad
\delta_n(A)=\dim_F\Delta_n(A).
\]
Thus $\delta_n(A)$ is the \emph{$n$th proper central codimension}.
The natural quotient map $P_n(A)\to P_n^z(A)$ has kernel
$\Delta_n(A)$, and hence
\begin{equation}\label{eq:codimension-relation}
c_n(A)=c_n^z(A)+\delta_n(A).
\end{equation}
In particular, $0\leq c_n^z(A)\leq c_n(A)$ for all $n\geq 1$.
If $A$ has no proper central polynomials, then
$c_n^z(A)=c_n(A)$ for every $n$.

We record some elementary consequences that will be used
throughout the paper. If $B\in\operatorname{var}(A)$ then $
\operatorname{Id}^z(A)\subseteq\operatorname{Id}^z(B),$ and therefore
$c_n(B)\leq c_n(A)$ and  $c_n^z(B)\leq c_n^z(A).$ Since these sequences depend only on the PI-equivalence class,
for $\mathcal V=\operatorname{var}(A)$ we may write
$c_n(\mathcal V)=c_n(A)$ and
$c_n^z(\mathcal V)=c_n^z(A)$.

Throughout this paper, we use the notation $\dotplus$ to denote a direct sum of vector spaces, whereas $\oplus$ denotes a direct sum of algebras.

\begin{remark} \label{directsum}
Let $A$ and $B$ be algebras. Evaluating polynomials componentwise and using $Z(A\oplus B)=Z(A)\oplus Z(B)$, we obtain
$\operatorname{Id}^z(A\oplus B)
=\operatorname{Id}^z(A)\cap\operatorname{Id}^z(B).$ Consequently, for every $n\geq 1$,
\[
\begin{aligned}
\max\{c_n^z(A),c_n^z(B)\}
&\leq c_n^z(A\oplus B)
\leq c_n^z(A)+c_n^z(B).
\end{aligned}
\]
Now, let $N$ be a nilpotent algebra with $N^t=\{0\}$. Hence
\[
c_n^z(N)=0
\quad\text{and}\quad
c_n^z(A\oplus N)=c_n^z(A)
\qquad\text{for all }n\geq t.
\]
Recall that a similar condition holds for ordinary codimensions.
\end{remark}

{\begin{remark}\label{rem:unital-nilpotent}
Let $R$ be a unital algebra such that
$R\sim_T D\oplus N$, for some algebra $D$ and a nilpotent algebra $N$.
Then $R\sim_T D$.
Indeed, choose $q\geq 1$ such that $N^q=0$.
For every multilinear polynomial
$f\in\operatorname{Id}(D)$, we have
\[
f(x_1,\ldots,x_m)x_{m+1}\cdots x_{m+q}
\in\operatorname{Id}(R).
\]
Evaluating the last $q$ variables at $1_R$
gives $f\in\operatorname{Id}(R)$.
Since the characteristic is zero, this yields
$\operatorname{Id}(D)\subseteq\operatorname{Id}(R)$.
The reverse inclusion follows from
$R\sim_T D\oplus N$.
\end{remark}}

The action of $S_n$ on $P_n$ by permutation of variables,
\[
\sigma\cdot f(x_1,\ldots,x_n)
=
f(x_{\sigma(1)},\ldots,x_{\sigma(n)}),
\]
preserves both $P_n\cap\operatorname{Id}(A)$ and
$P_n\cap\operatorname{Id}^z(A)$. Hence $P_n(A)$ and $P_n^z(A)$
are $FS_n$-modules. Their characters, denoted by $\chi_n(A)$ and
$\chi_n^z(A)$, are called the \emph{ordinary cocharacter} and the
\emph{central cocharacter} of $A$.

For a partition $\lambda\vdash n$, let $\chi_\lambda$ denote
the corresponding irreducible character of $S_n$ and let
$d_\lambda=\chi_\lambda(1)$. Complete reducibility gives
\[
\chi_n(A)
=
\sum_{\lambda\vdash n}m_\lambda(A)\chi_\lambda
\quad \mbox{ and }\quad 
\chi_n^z(A)
=
\sum_{\lambda\vdash n}m_\lambda^z(A)\chi_\lambda,
\]
where the multiplicities are nonnegative integers. Since
$P_n^z(A)$ is a quotient of $P_n(A)$, we have
$m_\lambda^z(A)\leq m_\lambda(A)$. Moreover,
\begin{equation*}\label{eq:central-character-degree} c_n(A)
=
\sum_{\lambda\vdash n}m_\lambda(A)d_\lambda \quad \mbox{ and }\quad 
c_n^z(A)
=
\sum_{\lambda\vdash n}m_\lambda^z(A)d_\lambda.
\end{equation*}

We shall also use highest weight vectors to study the central cocharacter sequence. The multiplicity $m_\lambda^z(A)$ equals the maximal number
of highest weight vectors of weight $\lambda$ that are
linearly independent modulo $\operatorname{Id}^z(A)$.
In particular, if a highest weight vector $f$ of weight
$\lambda$ is not central on $A$, then
\[
m_\lambda^z(A)>0
\qquad\text{and}\qquad
c_n^z(A)\geq d_\lambda.
\]
For the necessary background, we refer to \cite[Chapter 12]{Drensky}.

We say that $A$ has \emph{polynomial growth of central codimensions} if there exist $C>0$ and an integer $d\geq 0$
such that $c_n^z(A)\leq Cn^d$ for every $n\geq 1$. { This means that  $c_n^z(A)=\mathcal O(n^d)$.}
More precisely, we will use the notation $c_n^z(A)=\Theta(n^d)$ to mean that
there exist $C_1,C_2>0$ such that
\[
C_1n^d\leq c_n^z(A)\leq C_2n^d
\]
for all sufficiently large $n$. We say that $A$ has linear or quadratic central codimension growth when $c_n^z(A)=\Theta(n)$ or $c_n^z(A)=\Theta(n^2)$, respectively. 
We use the expression at most linear central codimension growth for the condition $c_n^z(A)=\mathcal{O}(n)$, and bounded central codimensions for the condition $c_n^z(A)=\mathcal O(1)$. We use analogous terminology for ordinary codimensions.

Motivated by the classical notion of minimal varieties
of polynomial codimension growth, we introduce an
analogous concept for central codimensions. 
For an integer $d\geq0$, a variety $\mathcal V$ is called $z$-minimal of degree $d$ if $c_n^z(\mathcal V)=\Theta(n^d)$ and every proper subvariety $\mathcal U$ of $\mathcal V$ satisfies either $c_n^z(\mathcal U)=0$ for all sufficiently large $n$, or $c_n^z(\mathcal U)=\Theta(n^t)$ for some integer
$0\leq t<d$.


\black
In particular, a variety of nonzero constant central codimension
growth is $z$-minimal if every proper subvariety has eventually zero
central codimensions.

A fundamental result of Regev \cite{RG} states that the codimension sequence of any PI-algebra is exponentially bounded. This means that, for a non-nilpotent algebra, there exist constants $b, C>0$ such that
\[1\leq c_n(A)\leq Cb^n
\qquad\text{for all }n\geq 1.
\]
This property motivates the study of the exponential rate of growth of the codimensions through the limit \[
\operatorname{exp}(A)
=
\lim_{n\to\infty}\sqrt[n]{c_n(A)},
\]
which, whenever it exists, is called the \emph{PI-exponent} of $A$.
Giambruno and Zaicev proved that, over a field of characteristic
zero, this limit exists for every PI-algebra and is a nonnegative
integer; see \cite{GiambrunoZaicev}.

In the central setting, \eqref{eq:codimension-relation} gives $0\leq c_n^z(A)\leq c_n(A)$, so the central codimension sequence
of a PI-algebra is also exponentially bounded. It is therefore
natural to consider the analogous limit
\[
\operatorname{exp}^z(A)
=
\lim_{n\to\infty}\sqrt[n]{c_n^z(A)},
\]
called the \emph{central exponent} of $A$, provided that it exists.
The existence and integrality of this limit, together with its
relation to the ordinary PI-exponent, were established by
Giambruno and Zaicev in \cite{GiZai2}, as follows.

\begin{theorem}[{\cite{GiZai2}}]
\label{thm:central-exponent}
Let $A$ be a PI-algebra over a field of characteristic zero.
Then $\operatorname{exp}^z(A)$ exists and is a nonnegative
integer. Moreover, $\operatorname{exp}^z(A)=\operatorname{exp}(A)$ if $\operatorname{exp}(A)\geq 2,$
whereas $\operatorname{exp}^z(A)\in\{0,1\}$
if $\operatorname{exp}(A)=1.$ If $\operatorname{exp}(A)=0$, then
$\operatorname{exp}^z(A)=0$.
\end{theorem}

Together with the polynomial-exponential dichotomy for ordinary codimensions, this theorem gives the following equivalence.

\begin{remark}
\label{rem:polynomial-growth-equivalence}
A PI-algebra has polynomial growth of ordinary codimensions
if and only if it has polynomial growth of central codimensions.
\end{remark}

Let $\mathcal G$ denote the infinite-dimensional Grassmann algebra over $F$, with identity element and generators
$\{e_1,e_2,\ldots\}$ satisfying
$e_ie_j=-e_je_i$ for all $i,j\geq 1.$
We also denote by $UT_2$ the algebra of $2\times 2$ upper
triangular matrices over $F$. These two algebras characterize
the boundary between polynomial and exponential growth of
ordinary codimensions. Combining the results of
\cite{GiamZai,Kem} with
Remark~\ref{rem:polynomial-growth-equivalence}, we obtain
the following central analogue.

{Notice that if $A$ is not a PI-algebra, then
$\operatorname{Id}^z(A)=\{0\}$. Indeed, a nonzero central
polynomial $f$ would yield a nonzero polynomial identity
$[f,x]$, where $x$ is a new variable.
Consequently, $c_n^z(A)=n!$ for every $n\geq 1$.
In particular, polynomial growth of central codimensions
implies that $A$ is a PI-algebra.}

\begin{proposition}\label{teokemer}
Let $A$ be an algebra over a field $F$ of characteristic zero.
The following conditions are equivalent:
\begin{enumerate}
    \item[(1)]
    $A$ has polynomial growth of central codimensions.

    \item[(2)]
    $\mathcal G,UT_2\notin\operatorname{var}(A)$.

    \item[(3)] $A\sim_T B_1\oplus\cdots\oplus B_m$, where $B_1,\ldots,B_m$ are finite-dimensional $F$-algebras such that $\dim_F(B_i/J(B_i))\leq 1$, for all
    $1\leq i\leq m,$ and $J(B_i)$ is the Jacobson radical of $B_i$.
\end{enumerate}
\end{proposition}

{ Each non-nilpotent algebra $B_i$ in Proposition~\ref{teokemer} admits
a Wedderburn--Malcev decomposition
\[
B_i=F\epsilon_i\dotplus J(B_i),
\]
where $\epsilon_i^2=\epsilon_i$ and $F\epsilon_i\cong F$.
We identify $F\epsilon_i$ with $F$ and write
$B_i=F\dotplus J(B_i)$ when no confusion can arise. Notice that $\epsilon_i$ is the identity element of $F\epsilon_i$, but is not necessarily an identity element of $B_i$. 
In our proofs, we will use 1 interchangeably to denote the element $\epsilon_i$ for algebras of type $F+J(B_i)$. }

As a consequence of the previous proposition, the central codimensions of a PI-algebra are
either polynomially bounded or grow exponentially; no intermediate growth occurs. This dichotomy motivates a finer
classification within the polynomial case. In the following
sections, we study varieties with bounded or linear central
codimensions and { determine the $z$-minimal varieties exhibiting
the corresponding growth thresholds.}


\section{Algebras with low-degree central codimension growth}

In this section, we introduce the algebras that will be used in the
classification results. We recall the relevant descriptions of their
polynomial identities and central polynomials, and determine the
central codimension growth needed in the subsequent arguments.

For an integer $m\geq 2$, let $UT_m$ denote the algebra of all
$m\times m$ upper triangular matrices over $F$. Let $e_{ij}$ be the
usual matrix units and set
\[
E=\sum_{i=1}^{m-1}e_{i,i+1},
\qquad
I_m=\sum_{i=1}^{m}e_{ii}.
\]
We consider the following subalgebras of $UT_m$:
\begin{align*}
&A_{m,1}
=\operatorname{span}_F
  \{e_{11},E,\ldots,E^{m-2},e_{12},e_{13},\ldots,e_{1m}\},\\
&A_{m,1}^*
=\operatorname{span}_F
  \{e_{mm},E,\ldots,E^{m-2},e_{1m},e_{2m},\ldots,e_{m-1,m}\},\\
&N_m
=\operatorname{span}_F
  \{I_m,E,\ldots,E^{m-2},e_{12},e_{13},\ldots,e_{1m}\}.
\end{align*}
When $m=2$, the list $E,\ldots,E^{m-2}$ is understood to be empty.

For $k\geq 1$, let
\[
\mathcal G_{2k}
=
\langle 1,e_1,\ldots,e_{2k}\mid e_ie_j=-e_je_i\rangle
\]
be the unital subalgebra of the infinite-dimensional Grassmann
algebra $\mathcal G$ generated by $1, e_1,\ldots,e_{2k}$.
In the particular case $k=1$, we have $\mathcal G_2\sim_T N_3$.

We define two families of algebras $\mathcal{F}_1=\{N, F, N_r,$ $ A_{m,1},$ $ A_{m,1}^*, UT_2
      \mid r\geq 3,\ m\geq 2\}$ and $\mathcal{F}_2=\{N,F,\mathcal G_{2k},\mathcal G\mid k\geq 1\}.$ Such families of algebras are essential in La Mattina's classification of the
subvarieties of $\operatorname{var}(\mathcal G)$ and
$\operatorname{var}(UT_2)$, which we recall for later use.

\begin{theorem}[{\cite[Theorems 5.2 and 5.4]{Danielasub}}]
\label{subvarieties}
Let $A$ be an algebra over a field $F$ of characteristic zero.
\begin{enumerate}
    \item If $A\in\operatorname{var}(UT_2)$ then $A$ is
    PI-equivalent to a finite direct sum of algebras in $\mathcal{F}_1.$
    \item If $A\in\operatorname{var}(\mathcal G)$ then $A$ is
    PI-equivalent to a finite direct sum of algebras in $\mathcal{F}_2.$
\end{enumerate}
Here $N$ denotes a nilpotent algebra.
\end{theorem}

\black

Next we introduce some algebras previously studied in \cite{GLa} and \cite{NascimentoCota}. First we consider the two-dimensional subalgebras of $UT_2$:
\[
A_1=A_{2,1}=Fe_{11}+Fe_{12}
\quad \mbox{ and }\quad 
A_1^*=A_{2,1}^*=Fe_{22}+Fe_{12}
\]
and notice that their centers are zero, so their ordinary and central
codimension sequences coincide.

\begin{lemma}\label{A1A1*} \cite[Lemma 3]{GLa}
For every $n\geq 1$, $c_n^z(A_1)=c_n^z(A_1^*)=n.$
\end{lemma}

Now we consider the following subalgebras of $UT_3$:
\begin{align*}
A_2&=N_3=F(e_{11}+e_{22}+e_{33})+Fe_{12}+Fe_{13}+Fe_{23},\\
A_4&=A_{3,1}=Fe_{11}+Fe_{12}+Fe_{13}+Fe_{23},\\
A_4^*&=A_{3,1}^*=Fe_{33}+Fe_{12}+Fe_{13}+Fe_{23},\\
A_5&=Fe_{22}+Fe_{12}+Fe_{13}+Fe_{23},\\
A_6&=F(e_{11}+e_{33})+Fe_{12}+Fe_{13}+Fe_{23}.
\end{align*}

The next lemma summarizes the information on central polynomials 
needed for these algebras.

\begin{lemma}\label{lem:UT3-central}
The following statements hold.
\begin{enumerate}
    \item The algebras $A_4$ and $A_4^*$ have no proper central
    polynomials. Consequently,
$ \operatorname{Id}^z(A_4)=\operatorname{Id}(A_4)
      =\langle [x_1,x_2]x_3x_4\rangle_T$ and $
    \operatorname{Id}^z(A_4^*)
    =\operatorname{Id}(A_4^*)
      =\langle x_3x_4[x_1,x_2]\rangle_T.$
    \item For every $n\geq 4$, $c_n^z(A_4)=c_n^z(A_4^*)=n(n-1)$  and $c_n^z(A_5)=c_n^z(A_6)=2n-1.$ \end{enumerate}
\end{lemma}

\begin{proof}
The centers of $A_4$ and $A_4^*$ are zero. Thus the first two
statements follow from the descriptions of their ordinary identities and codimensions in \cite[Lemma 6]{GLa}.
The last statement is the central codimension formula obtained in \cite[Lemma 5.3]{NascimentoCota}.
\end{proof}

The generators of $\operatorname{Id}(A_6)$ will be used explicitly in the structural arguments. By \cite[Theorem 3.1]{Sandra},  $\operatorname{Id}(A_6)=\langle f_1,\ldots,f_6\rangle_T,$ where
\begin{align*}
f_1 \; &=\; { St_4(x_1,x_2,x_3,x_4)= }[x_1,x_2][x_3,x_4] + [x_3,x_4][x_1,x_2] + [x_2,x_3][x_1,x_4] + [x_1,x_4][x_2,x_3] \\ 
& \quad + [x_3,x_1][x_2,x_4] + [x_2,x_4][x_3,x_1], \\[6pt]
f_2 \; &=\; [x_1,x_2][x_3,x_4] + [x_3,x_4][x_1,x_2] + x_1[x_3,x_4]x_2 - x_2[x_3,x_4]x_1, \\[6pt]
f_3 \; &=\; 2x_1x_2[x_4,x_3] + x_4x_1[x_3,x_2] + x_4x_2[x_3,x_1] + x_3x_1[x_2,x_4] + x_3x_2[x_1,x_4] \\
&\quad + [x_1,x_2][x_3,x_4] + [x_1,x_3][x_4,x_2] - [x_4,x_2][x_1,x_3]
+ [x_2,x_3][x_4,x_1] - [x_4,x_1][x_2,x_3], \\[6pt]
f_4 \; &=\; [x_1,x_2]x_3[x_4,x_5],\quad  f_5 \; =\; [[x_1,x_2][x_3,x_4],x_5] \quad \mbox{and}\quad f_6 \; =\; [x_3[x_1,x_2]x_4,x_5]. \
\end{align*}

We shall also need the following subalgebras of $UT_4$:
\begin{align*}
A_7={}&F(e_{11}+e_{22}+e_{33}+e_{44})+Fe_{12}+Fe_{13}+Fe_{14}+Fe_{23}+Fe_{24}+Fe_{34},\\[4pt]
A_8={}&F(e_{11}+e_{22}+e_{33})+Fe_{12}+Fe_{13}+Fe_{14}+Fe_{23}+Fe_{24}+Fe_{34},\\[4pt]
A_8^*={}&F(e_{22}+e_{33}+e_{44})+Fe_{12}+Fe_{13}+Fe_{14}+Fe_{23}+Fe_{24}+Fe_{34}.
\end{align*}

The algebras $A_7$ and $N_4$ are PI-equivalent
\cite[Remark 3.7]{NascimentoCota}. Consequently, their spaces
of central polynomials coincide, and the known results for
$N_4$ \cite[Theorem 3.3]{AnaThaisWillerJuan} give the following description. Moreover, $A_8$ and $A_8^*$ have trivial center and so the corresponding result comes from \cite[Lemma 9]{Daniela}.

\begin{lemma} \label{centralcodimensionA7} The following statements hold.
    \begin{enumerate}
        \item $\textnormal{Id}^z(A_7)=\textnormal{Id}^z(N_4)=\textnormal{Id}(N_3)=\langle [x_1,x_2,x_3],[x_1,x_2][x_3,x_4]\rangle_T$. 
        \item $\textnormal{Id}^z(A_8)=\langle [x_1,x_2][x_3,x_4]x_5,[x_1,x_2,x_3]x_4\rangle_T$ and $\textnormal{Id}^z(A_8^*)=\langle x_1[x_2,x_3][x_4,x_5],x_1[x_2,x_3,x_4]\rangle_T$.
        
        \item $c_n^z(A_7)=c_n^z(N_4)=c_n(N_3)= 1+\binom{n}{2}$ and $c_n^z(A_8)=c_n^z(A_8^*)= n+\binom{n}{2}(n-2)$, for all $n> 4$.
    \end{enumerate}
\end{lemma}
The finite-dimensional Grassmann algebras provide another family
whose central codimensions have polynomial growth. We recall both
their central polynomials and the corresponding codimension formula, see \cite[Proposition 6]{PKS} and \cite[Theorem 4.4]{AnaThaisWillerJuan}.

\begin{lemma}\label{G2k}
For every $k\geq 1$, $\operatorname{Id}^z(\mathcal G_{2k})
=\bigl\langle [x_1,x_2],\;y[x_1,x_2,x_3],y[x_1,x_2][x_3,x_4]\cdots[x_{2k-1},x_{2k}]
\bigr\rangle^T.$ Moreover, for every $n\geq 1$,
\[
c_n^z(\mathcal G_{2k})
=\sum_{i=0}^{k-1}\binom{n-1}{2i}.
\]
\end{lemma}

In particular, $c_n^z(A_2)=c_n^z(\mathcal G_2)=1$ and $c_n^z(\mathcal G_4)=1+\binom{n-1}{2},$ for all $ n\geq 1.$
Thus $A_2$ has bounded central codimensions, whereas $\mathcal G_4$ has quadratic central codimension growth.

Finally, consider

\[
K=
\left\{
\begin{pmatrix}
a&b&c&d\\
0&0&e&f\\
0&0&a&g\\
0&0&0&a
\end{pmatrix}
\mid  a,\ldots ,g\in F
\right\}.
\]

Writing $\epsilon_K=e_{11}+e_{33}+e_{44}$, we have $K=F\epsilon_K\dotplus J(K),$ where
$J(K)=\operatorname{span}_F\{e_{ij}\mid 1\leq i<j\leq4\}.$ Direct multiplication gives
$Z(K)=Fe_{14}=J(K)^3.$

Now consider a finite-dimensional algebra  $A=F\dotplus J$ with Jacobson radical $J$.
Let $q$ denote the nilpotency index of $J$, so that $J^q=0$.
In particular, $J^q\subseteq Z(A)$, and hence there always exists
an integer $k\geq 1$ such that $J^k\subseteq Z(A)$.
Since $A/J\cong F$, the standard estimate for ordinary codimensions
gives $c_n(A)=O(n^{q-1})$, and consequently
$c_n^z(A)=O(n^{q-1})$.

The following lemma provides a refinement of this bound.
Indeed, if $J^k\subseteq Z(A)$ for some $k<q$, then the degree
of the polynomial bound for the central codimensions can be
reduced from $q-1$ to $k-1$. Thus, the result takes into account
not only the nilpotency of the radical, but also its interaction
with the center of $A$.

\begin{lemma}\label{lem:central-radical-bound}
Let $A=F\dotplus J$ be a finite-dimensional algebra over a field $F$
of characteristic zero, where $J=J(A)$.
If $J^k\subseteq Z(A)$ for some integer $k\geq 1$, then there
exists a constant $\alpha>0$ such that
\[
  c_n^z(A)\leq \alpha n^{k-1}
\]
for all $n\geq 1$.
\end{lemma}

\begin{proof}
Set $B=A/J^k$. Every polynomial identity of $B$ takes its values
in $J^k\subseteq Z(A)$ when evaluated on $A$. Hence $\operatorname{Id}(B)\subseteq\operatorname{Id}^z(A)$ and $c_n^z(A)\leq c_n(B).$

Write $B=F+ R$, where $R=J/J^k$ and $R^k=0$,
and fix a basis $\{1,r_1,\ldots,r_d\}$ of $B$,
with $r_1,\ldots,r_d\in R$.
Every evaluation of a multilinear polynomial on this basis
involving at least $k$ radical elements vanishes.
Set $t=\min\{n,k-1\}$. There are
\[
T_n=\sum_{r=0}^{t}\binom{n}{r}d^r
\]
basis substitutions involving fewer than $k$ radical
elements. The corresponding evaluation maps define a
linear map $P_n\rightarrow B^{T_n}$ whose kernel is $P_n\cap\operatorname{Id}(B)$,
by multilinearity.

Consequently, for every $n\geq1$,
\[
\begin{aligned}
c_n^z(A)\leq c_n(B)\leq (\dim_F B)
       \sum_{r=0}^{t}\binom{n}{r}d^r
\leq
(\dim_F B)
\left(\sum_{r=0}^{k-1}\frac{d^r}{r!}\right)n^{k-1}.
\end{aligned}
\]
\end{proof}

\begin{lemma} \label{k}
For the algebra $K$ we have $c_n^z(K)=\Theta(n^2)$.
\end{lemma}

\begin{proof}Since $J^3=Fe_{14}=Z(K),$ Lemma~\ref{lem:central-radical-bound}, applied with $k=3$, gives $c_n^z(K)\leq\alpha_2n^2$ for some $\alpha_2>0$.

In order to establish the lower bound, for all $n\geq 4$, consider the highest weight vector $f=[x_1,x_2]^2x_1^{n-4}$ associated to the {partition} $(n-2,2)$. By evaluating $x_{1}\mapsto e_{11}+e_{33}+e_{44}$ and $x_2\mapsto e_{12}+e_{23}$ we obtain $f\mapsto -e_{13}$ which is a noncentral element. Therefore, $$\frac{n(n-3)}{2}=\deg \chi_{(n-2,2)}\leq c_n^z(K),$$ and the proof follows.     
\end{proof}

The preceding results exhibit algebras with several small central codimension growth rates. The algebra $\mathcal G_2$ has bounded central codimensions; $A_1$, $A_1^*$, $A_5$ and $A_6$ have linear growth;
$A_4$, $A_4^*$, $A_7$, $\mathcal G_4$ and $K$ have quadratic growth; and $A_8$ and $A_8^*$ have cubic growth. These computations will be used in the classification arguments of the next sections.

\section{When the central codimensions are bounded by a constant}

In this section, we classify algebras whose central codimension
sequences are bounded by a constant. By Proposition~\ref{teokemer}, it suffices
to study finite direct sums of nilpotent algebras and
finite-dimensional algebras with a one-dimensional semisimple component.

Let $A=F\dotplus J(A)$ be a finite-dimensional algebra. Here we identify $F$ with a fixed
one-dimensional semisimple subalgebra of $A$, and write $1_F$ for
its identity element. Notice that $1_F$ need not be an identity
element of $A$.

Set $J=J(A)$. By \cite[Lemma 9.7.1]{GiambrunoZai}, the radical admits the decomposition into vector subspaces
\begin{equation} \label{decompj}
    J(A)=J_{10}\dotplus J_{01} \dotplus J_{00}\dotplus  J_{11},
\end{equation}
where $J_{ik}
=
\{a\in J\mid 1_Fa=ia,\ a1_F=ka\},$ $ i,k\in\{0,1\}.$ Moreover, $J_{ik}J_{rs}\subseteq\delta_{kr}J_{is}$ where $\delta_{kr}$ is the Kronecker delta function.

We next recall the results established in
\cite[Lemmas~9, 10, 15, 20, 21]{GLa} and in \cite[Lemmas~13 and~16]{Daniela} concerning finite-dimensional algebras of the form $A=F\dotplus J$, where $J=J(A)$ has a decomposition as given in (\ref{decompj}).

\begin{lemma} \label{1} For the algebra $A=F\dotplus J$ we have the following.

\begin{enumerate}
 \item[(1)]If  $A_{1}\notin \var(A)$ then $J_{10}= \{0\}$.

  \item[(2)]If $A_{1}^*\notin \var(A)$ then $J_{01}= \{0\}$.

   \item[(3)]If $A_{2}\notin \var(A)$ then $[J_{11}, J_{11}] = \{0\}$.

   \item[(4)] If $A_4\notin \var(A)$ then $J_{10}J_{00}= \{0\}$.
    \item[(5)] If $A_4^*\notin \var(A)$ then $J_{00}J_{01}=\{0\}$.
\item[(6)]  If $J_{10}J_{00}=J_{00}J_{01}= \{0\}$ and $J_{01}J_{10}\neq \{0\}$ then $A_5\in \var(A).$

\item[(7)] If $J_{10}J_{01}\neq \{0\}$ and $J_{01}J_{10}=\{0\}$ then $A_6\in \var(A).$

 \item[(8)]     If $[J_{11},J_{11},J_{11}]\neq \{0\}$ then $A_7\in \var(A)$.

        \item[(9)] If $\mathcal{G}_4, A_7, A_8,A_8^*\notin \var(A)$ and $J_{00}=J_{10}J_{01}=\{0\}$ then $A\in \var(UT_2).$
    \end{enumerate}
\end{lemma}

\black

We can now state the first
characterization of varieties whose sequence of central codimensions is bounded by a constant.

\begin{theorem} \label{boundedbyconstant}
Let $A$ be an algebra over a field $F$ of characteristic zero.
Then the following conditions are equivalent:
\begin{enumerate}
    \item There exists a constant $\alpha >0$ such that $c_n^z(A)\leq \alpha$, for all $n\geq 1.$

    \item $A$ is PI-equivalent to one of the algebras  $N$, $F\oplus N$ or $\mathcal G_2\oplus N,$ where $N$ is a nilpotent $F$-algebra, possibly zero.

    \item $A_1,A_1^*,\mathcal G_4,A_7
    \notin\operatorname{var}(A).$
\end{enumerate}
\end{theorem}

\begin{proof}
First we notice that Lemmas~\ref{A1A1*}, \ref{centralcodimensionA7} and \ref{G2k} imply that $A_1, A_1^*, \mathcal{G}_{4}$ and $A_7$ have at least linear growth of central codimensions. Consequently, by the monotonicity of central codimensions
under inclusion of varieties, we obtain $(1)\Rightarrow(3)$.

Suppose item (2) is true and choose $s\geq 1$ such that $N^s=\{0\}$.
For all $n\geq s$, we have
\[
c_n^z(N)=c_n^z(F\oplus N)=0
\qquad\text{and}\qquad
c_n^z(\mathcal G_2\oplus N)=c_n^z(\mathcal G_2)=1.
\]
These sequences are therefore eventually constant and hence bounded.
This proves $(2)\Rightarrow(1)$, so it remains only to establish
$(3)\Rightarrow(2)$.

Assume that $A_1, A_1^*, \mathcal{G}_{4}, A_7 
\notin \operatorname{var}(A).$ Since $\mathcal{G}_{4}\in \operatorname{var}(\mathcal{G})$  and $A_1\in \operatorname{var}(UT_2),$ it follows that $\mathcal{G},UT_2\notin \operatorname{var}(A).$
Therefore, by Proposition~\ref{teokemer},  $A$ has polynomial growth of central codimensions. Moreover, after replacing $A$ by a PI-equivalent algebra, we may assume that
\[
A=B_1\oplus\cdots\oplus B_m,
\]
where each $B_i$ is a finite-dimensional algebra over $F$ such that either $B_i$ is nilpotent or $B_i=F\dotplus J(B_i),$ where $J(B_i)$ denotes the Jacobson radical of $B_i$.

If every $B_i$ is nilpotent, then $A$ is nilpotent, since it is a finite direct sum of nilpotent algebras. Hence, in this case, there is nothing more to prove.
So, we may assume that some $B_i$ is a non-nilpotent summand and consider the decomposition
\[
J(B_i)
 =
J_{10}^{(i)}\dotplus  J_{01}^{(i)}
\dotplus J_{00}^{(i)}\dotplus  J_{11}^{(i)}
\]
given in \eqref{decompj}. Since $B_i\in\operatorname{var}(A)$ and
$A_1,A_1^*\notin\operatorname{var}(A),$
we also have $A_1,A_1^*\notin\operatorname{var}(B_i).$ Thus, by items (1) and (2) of Lemma~\ref{1}, we have $J_{10}^{(i)}=J_{01}^{(i)}=\{0\}.$ Consequently, $B_i=Q_i\oplus N_i,$ where $Q_i=F\dotplus J_{11}^{(i)}$ and $N_i=J_{00}^{(i)}$ is a nilpotent algebra.

Since $Q_i\in\operatorname{var}(A)$, we have
$A_7\notin\operatorname{var}(Q_i)$ and Lemma~\ref{1} implies that $[x_1,x_2,x_3]\equiv 0$ on $Q_i$. Hence, $Q_i\in\operatorname{var}(\mathcal{G}).$ Furthermore, since $\mathcal{G}_4\notin\operatorname{var}(A)$, we also have $\mathcal{G}_4\notin\operatorname{var}(Q_i).$
It now follows from  Theorem \ref{subvarieties} that $Q_i$ is PI-equivalent to either $F\oplus N_i'$ or $\mathcal{G}_2\oplus N_i',$ where $N_i'$ is nilpotent.

Since $B_i=Q_i\oplus N_i$, we conclude that every non-nilpotent summand $B_i$ is PI-equivalent to either $F\oplus(N_i'\oplus N_i)$ or $\mathcal{G}_2\oplus(N_i'\oplus N_i).$ In each case, $N_i'\oplus N_i$ is nilpotent.

Finally, applying this argument to every non-nilpotent summand of $A$ and collecting all the nilpotent summands, we obtain a nilpotent algebra $N$ such that one of the following situations occurs:
$
A\sim_T N,
A\sim_T F\oplus N,
A\sim_T\mathcal{G}_2\oplus N.
$
This proves $(2)$ and completes the proof.
\end{proof}

Two algebras $A$ and $B$ are said to be \emph{comparable}
if $\operatorname{Id}(A)
\subseteq \operatorname{Id}(B)$ or $\operatorname{Id}(B)
\subseteq \operatorname{Id}(A).$
They are said to be \emph{incomparable} if they are not comparable. As an example, the polynomials $[x_1,x_2]x_3$ and
$x_3[x_1,x_2]$ show that the algebras $A_1$ and $A_1^*$ are incomparable.

\begin{corollary}
The varieties $\operatorname{var}(A_1)$ and $\operatorname{var}(A_1^*)$ are  $z$-minimal of degree 1. 
\end{corollary}
\begin{proof}
     Let $\mathcal U$ be a proper subvariety of either $\operatorname{var}(A_1)$ or $\operatorname{var}(A_1^*)$. Since $A_1$ and $A_1^*$ are incomparable, $\mathcal U$ contains neither $A_1$ nor $A_1^*$. Moreover, since $c_n^z(\mathcal U)$  is bounded by a linear function, $\mathcal U$ also contains neither $\mathcal G_4$ nor $A_7$, both of which have quadratic central codimension growth. Therefore, Theorem~\ref{boundedbyconstant} implies that the central codimensions of $\mathcal U$ are bounded by a constant.
\end{proof}

\black
\begin{corollary} 
The variety $\operatorname{var}(\mathcal G_2)$ is the unique
$z$-minimal variety of nonzero constant central codimension growth.
The only $z$-minimal varieties of degree $1$
are $\operatorname{var}(A_1)$ and $\operatorname{var}(A_1^*)$.
\end{corollary}

\begin{proof}
By Theorem~\ref{boundedbyconstant}, every variety with bounded central codimensions
that are not eventually zero is generated by
$\mathcal G_2\oplus N$, for some nilpotent algebra $N$.
The $z$-minimality therefore forces it to equal
$\operatorname{var}(\mathcal G_2)$.
On the other hand, the same theorem implies that every proper subvariety
of $\operatorname{var}(\mathcal G_2)$ has eventually zero central
codimensions.

Now let $\mathcal V$ be a $z$-minimal variety of linear central codimension growth. Since $\mathcal G_4$ and $A_7$ have quadratic central codimension growth, neither belongs to $\mathcal V$.
If neither $A_1$ nor $A_1^*$ belonged to $\mathcal V$, Theorem~\ref{boundedbyconstant} would imply that its central codimensions are bounded, a contradiction. Thus, $A_1\in\mathcal V$
or $A_1^*\in\mathcal V$. The $z$-minimality therefore forces $\mathcal V=\operatorname{var}(A_1)$ or
$\mathcal V=\operatorname{var}(A_1^*)$.
\end{proof}

\section{Varieties with at most linear growth of central codimensions}

In this section, we classify the varieties whose central codimensions are bounded by a linear function. More precisely, we show that every such variety is generated by either a nilpotent algebra or the direct sum of a nilpotent algebra and an algebra from a fixed finite family. We also determine the { $z$-minimal varieties of degree 2.}
To achieve our goals, we consider $$\mathcal{I}=\{\mathcal G_{4},A_{7},A_{4},A_{4}^{*}, A_8, A_8^*\}.$$

In the next three lemmas, we will consider 
$A=F\dotplus J$ to be a finite-dimensional algebra over a field $F$
of characteristic zero, where $J=J(A)$. 

\black
{\begin{lemma}\label{lemma-A5}
Suppose that $J_{10}J_{01}=\{0\}$ and $J_{01}J_{10}\neq \{0\}$. If $Q\notin\textnormal{var}(A)$, for all $Q\in \mathcal{I}$, then $A$ is PI-equivalent
to either $A_{5}\oplus N$ or $
A_{5}\oplus \mathcal{G}_2\oplus N,$ where $N$ is a nilpotent algebra. 
\end{lemma}

\begin{proof}
Consider the subalgebras $B=F\dotplus J_{11}$ and $
D=F\dotplus J_{10}\dotplus J_{01}\dotplus J_{11}\dotplus I$ of $A$, 
where $\{0\}\neq I=J_{01}J_{10}\subseteq J_{00}$, by hypothesis. Since $A_{4},A_{4}^{*}\notin\operatorname{var}(A),$ by Lemma~\ref{1}, we get $J_{00}J_{01}=J_{10}J_{00}=\{0\}$. 
On the other hand, since $J_{10}J_{01}=\{0\}$ and  $J_{01}J_{10}\neq\{0\},$ we use again Lemma~\ref{1}  to conclude that $A_{5}\in\operatorname{var}(D).$ Although $D\cap J_{00}=I$ may be nonzero, we have
$A=D+J_{00}$ and $DJ_{00}=J_{00}D=\{0\}$.
Thus, by multilinearity, a multilinear polynomial is an identity of
$A$ if and only if it is an identity of both $D$ and $J_{00}$.
Since the characteristic is zero, it follows that
\[
\operatorname{Id}(A)
=
\operatorname{Id}(D)\cap\operatorname{Id}(J_{00})
=
\operatorname{Id}(D\oplus J_{00}).
\]
Consequently, $A\sim_T D\oplus J_{00}$.

Since $A_{7}\notin\operatorname{var}(A)$, by Lemma \ref{1}, we have $[J_{11},J_{11},J_{11}]=\{0\}.$
{ Thus the unital algebra $B=F\dotplus J_{11}$ satisfies
$[x_1,x_2,x_3]\equiv 0$. The exclusion
$\mathcal G_4\notin\operatorname{var}(B)$ and
Theorem~\ref{subvarieties} imply that
$B\sim_T F\oplus N'$ or
$B\sim_T\mathcal G_2\oplus N'$,
where $N'$ is nilpotent.
By Remark~\ref{rem:unital-nilpotent},
we conclude that
$B\sim_T F$ or $B\sim_T\mathcal G_2$.
In particular,
$[J_{11},J_{11}][J_{11},J_{11}]=\{0\}$.} 

If $B$ is commutative then 
$[J_{11},J_{11}]=\{0\}$ and in this case, 
using that $J_{10}J_{01}=\{0\}$, one immediately
obtains
$
[D,D]\subseteq J_{10}+J_{01}+I.
$
Furthermore, $[D,D]D\subset J_{01}+I$ and $D(J_{01}+I)=\{0\}$  and so
$D[D,D]D=\{0\}.$

Thus, $x_{1}[x_{2},x_{3}]x_{4}\equiv 0$
is a polynomial identity of $D$. By \cite[ Lemma~6]{GLa}, $\operatorname{Id}(A_5)
=
\langle x_1[x_2,x_3]x_4\rangle_T.$
Hence $D\in\operatorname{var}(A_5)$.
Together with $A_5\in\operatorname{var}(D)$, this gives
$D\sim_T A_5$. {Recalling that $A\sim_T D\oplus J_{00}$, we conclude in this case that
$A\sim_T A_5\oplus N$, with $N=J_{00}$ nilpotent.}

Now consider $B \sim_{T} \mathcal{G}_2$.  Since $A_8,A_8^*\notin\operatorname{var}(A)$,
it follows from \cite[Lemma~14]{Daniela} that
\begin{equation*}\label{eq-interactions-C}
[J_{11},J_{11}]J_{10}=\{0\}
\qquad\text{and}\qquad
J_{01}[J_{11},J_{11}]=\{0\}.
\end{equation*}

Let $H = [J_{11}, J_{11}] + [J_{11}, J_{11}] J_{11}$ and $L=J_{10}+J_{01}+I$. Notice that both $H$ and $L$ are  two-sided ideals of  $D$ and since $H\subseteq J_{11}$ and $L\subseteq J_{10} +J_{01} + J_{00}$, the direct decomposition of the radical gives $H\cap L=\{0\}.$ Therefore, the canonical homomorphism
\[
\Phi:D\rightarrow D/H\oplus D/L,
\qquad
\Phi(d)=(d+H,d+L),
\]
is injective. Both $D/H$ and $D/L$ are homomorphic images of $D$, whereas $D$ is
isomorphic to a subalgebra of their direct sum. Consequently,
\[
D\sim_{T}D/H\oplus D/L.
\]

In $D/H$, the image of $J_{11}$ is commutative, since $[J_{11},J_{11}]\subseteq H.$ Moreover, $I=J_{01}J_{10}$ survives in the quotient because $H\cap I=\{0\}.$ Thus, by the argument used in the first case, $D/H\sim_{T}A_{5}.$

Moreover, $D/L\cong F\dotplus J_{11}$ and since $B= F\dotplus J_{11}\sim_{T}\mathcal{G}_{2}$, it follows that
\[
D
 \sim_{T}D/H\oplus D/L
 \sim_{T}A_{5}\oplus \mathcal{G}_{2}.
\]

Finally, recalling that $A\sim_T D\oplus J_{00}$, we obtain
$A\sim_T A_5\oplus \mathcal G_2\oplus N$,
where $N=J_{00}$ is nilpotent.

\end{proof}}

\black


\begin{lemma}\label{commutativityj10j01}
Suppose that $J_{01}J_{10}=\{0\}$ and $[J_{10}J_{01},J_{11}]\neq\{0\}.$ Then $K\in\operatorname{var}(A)$.
\end{lemma}

\begin{proof} 
Choose \(a\in J_{10}\), \(b\in J_{01}\) and \(c\in J_{11}\)
such that \([ab,c]\ne0\). In particular, \(ab\ne0\).
Put \(t=1_F+c\in B=F+J_{11}\).
By the multiplication rules for the radical components,
\begin{equation}\label{equation1}
 1_Fa=a,\quad a1_F=0,\qquad 1_Fb=0,\quad b1_F=b,\qquad
 1_Ft=t1_F=t,\quad at=tb=ba=0.
\end{equation}
Moreover, since \(c\) is nilpotent, \(t\) is invertible in the unitary algebra \(B\). In particular, \(bt\ne0\).

We first observe that
\begin{equation} \label{equation2}
  abt\quad\text{and}\quad tab
  \quad\text{are linearly independent over }F.
\end{equation}

Indeed, choose \(k\ge2\) such that
\(ab\in J^k\setminus J^{k+1}\) and suppose that
\(
  \alpha abt+\beta tab=0, \alpha,\beta\in F.
\)

Since \(c\in J\), reducing modulo \(J^{k+1}\) yields
\((\alpha+\beta)ab\equiv0\pmod{J^{k+1}}\).
Thus \(\beta=-\alpha\), and consequently
\(
  0=\alpha(abt-tab)=\alpha[ab,c].
\)
The choice of \(a,b,c\) implies \(\alpha=\beta=0\), proving (\ref{equation2}).

Now considering the following generators of \(K\):
\(
 e_K=e_{11}+e_{33}+e_{44},\
 u=e_{12},\ v=e_{23},\ w=e_{34},
\)
we notice that \(  uv=e_{13}, \ vw=e_{24}, \ uvw=e_{14}\), so we have
\(
 K=\mbox{span}\{e_K,u,v,w,uv,vw,uvw\}.
\)
Next we will prove that \(\operatorname{Id}(A)\subseteq \operatorname{Id}(K)\).

It is enough to check that every multilinear identity of \(A\)
vanishes when its variables are evaluated in the generating set
\(\{e_K,u,v,w\}\).
For completeness, recall why this suffices:
if \(f\) is a multilinear identity of \(A\) and its variables are
replaced by words in these generators, replace each occurrence
of a generator in those words by a new, distinct variable.
The resulting polynomial is again a multilinear identity of \(A\),
since \(\Id A\) is a \(T\)-ideal.
Its evaluation at the original generators therefore vanishes.
Thus \(f\) vanishes on arbitrary words, and hence, by multilinearity,
on all of \(K\). As \(\operatorname{char}F=0\), this proves the desired
inclusion of \(T\)-ideals.

Let \(f=f(x_1,\ldots,x_n)\in \operatorname{Id}(A)\) be multilinear.
Consider an evaluation of \(f\) in which every variable takes one of
the values \(e_K,u,v,w\).
If any of \(u,v,w\) occurs more than once, every monomial evaluates
to zero: a nonzero product of the matrix units
$u, v, w$, with arbitrary
factors \(e_K\) inserted, cannot repeat any of these matrix units.
Therefore, we may assume that each of \(u,v,w\) occurs at most once.

 \begin{itemize}
    \item \textbf{Case 1.} No variable is evaluated at \(w\).
\end{itemize} 

The subalgebra of \(K\) generated by \(e_K,u,v\) has basis
\(\{e_K,u,v,uv\}\). By (1) and \(ab\ne0\), the assignment
\[
  e_K\longmapsto 1_F,\qquad
  u\longmapsto a,\qquad
  v\longmapsto b,\qquad
  uv\longmapsto ab
\]
defines an isomorphism onto the subalgebra of \(A\) generated
by \(1_F,a,b\). The elements \(1_F,a,b,ab\) are linearly independent, since they belong, respectively, to the distinct direct summands
\(F1_F,J_{10},J_{01}\), and \(J_{11}\) of \(A\).
Hence \(f\) vanishes under the evaluation under consideration.

 \begin{itemize}
    \item \textbf{Case 2.} \(w\) occurs, but neither \(u\) nor \(v\) occurs.
\end{itemize}

Since \(e_Kw=we_K=w\), every monomial of \(f\)
has value \(w\) under this evaluation.
On the other hand, evaluating all variables at \(1_F\in A\)
shows that the sum of the coefficients of \(f\) is zero.
Therefore the evaluation in \(K\) vanishes.

 \begin{itemize}
    \item \textbf{Case 3.} \(u\) and \(w\) occur, but \(v\) does not.
\end{itemize}

Every monomial vanishes in \(K\), because
\(uw=wu=0\), even when factors \(e_K\) are inserted.

 \begin{itemize}
    \item \textbf{Case 4.} \(v\) and \(w\) occur, but \(u\) does not.
\end{itemize}

A monomial has nonzero value in \(K\) precisely when
\(v\) is its first factor and \(w\) occurs later.
All such monomials have the same value \(vw\).
Let \(\gamma\) be the sum of their coefficients.
Thus the evaluation of \(f\) in \(K\) is \(\gamma vw\).

Now evaluate \(f\) in \(A\), replacing \(v\) by \(b\),
\(w\) by \(t\), and all occurrences of \(e_K\) by \(1_F\).
By (1), exactly the same monomials can have nonzero values,
each equal to \(bt\).
Hence \(0=\gamma bt\). Since \(bt\ne0\), we obtain
\(\gamma=0\), as required.

 \begin{itemize}
    \item \textbf{Case 5.} \(u\), \(v\) and \(w\) all occur.
\end{itemize}

In \(K\), a monomial has nonzero value precisely when
\(u\) is immediately followed by \(v\), and \(w\) occurs
after \(v\). All such monomials have value \(uvw\).
Let \(\alpha\) be the sum of their coefficients.
The evaluation of \(f\) in \(K\) is therefore
\begin{equation} \label{equation3}
  \alpha\,uvw.    
\end{equation}

Evaluate the same multilinear polynomial in \(A\), replacing
\(u,v,w,e_K\) by \(a,b,t,1_F\), respectively.
The multiplication rules (\ref{equation1}) imply that a monomial can have
nonzero value only when \(a\) is immediately followed by \(b\).
There are exactly two possible positions for \(t\):
it occurs after \(b\), giving the value \(abt\), or it occurs
before \(a\), giving the value \(tab\).
All other factors are copies of \(1_F\) and do not change these values.
Consequently,
\[
  0=f(\ldots,a,\ldots,b,\ldots,t,\ldots)
   =\alpha\,abt+\beta\,tab,
\]
where \(\beta\in F\) is the sum of the coefficients of the
monomials in which \(t\) occurs before the adjacent pair \(ab\).
By the linear independence established in (\ref{equation2}), we obtain
\(\alpha=\beta=0\).
Equation (\ref{equation3}) shows that the evaluation of \(f\) in \(K\) vanishes.

\smallskip
These cases exhaust all evaluations in the generating set and so, $f\in \operatorname{Id}(K)$.
We conclude that \(K\in\var(A)\).

\end{proof} 

Let us consider the set $$\mathcal{J}=\mathcal{I}\cup \{K \}= \{\mathcal G_{4},A_{7},A_{4},A_{4}^{*}, A_8, A_8^*, K\}.$$

\begin{lemma}\label{lemma-A6} Suppose
$J_{00}J_{01}=J_{10}J_{00}=J_{01}J_{10}=\{0\}.$ If $J_{10}J_{01}\neq \{0\}$ and $Q\notin\textnormal{var}(A)$, for all $Q\in \mathcal{J}$, then $A$ is PI-equivalent
to either $A_{6}\oplus N$ or { $
A_{6}\oplus \mathcal{G}_2\oplus N,$} where $N$ is a nilpotent algebra. 
\end{lemma}
\begin{proof}
Set $R=F\dotplus J_{11}$ and $C=[J_{11},J_{11}].$ Since $A_{7}\notin\operatorname{var}(A)$,  by Lemma \ref{1}, $R$ satisfies the identity $[x_{1},x_{2},x_{3}]\equiv 0.$
{Moreover, since
$\mathcal G_4\notin\operatorname{var}(A)$,
Theorem~\ref{subvarieties} implies that
$R\sim_T F\oplus N'$ or
$R\sim_T \mathcal{G}_2\oplus N'$,
where $N'$ is nilpotent.
Since $R$ is unital,
Remark~\ref{rem:unital-nilpotent} gives
$R\sim_T F$ or $R\sim_T \mathcal{G}_2$.}
In particular,
\[
C^2=[J_{11},J_{11}][J_{11},J_{11}]
=\{0\}.
\]

Since $A_8,A_8^*\notin \textnormal{var}(A)$, by \cite[Lemma~14]{Daniela}, we therefore obtain 
\begin{equation}\label{eq:annihilation-A6}
CJ_{10}=\{0\}
\qquad\text{and}\qquad
J_{01}C=\{0\}.
\end{equation}

Consider the subalgebra
\[
D=F\dotplus J_{10}\dotplus J_{01}\dotplus J_{11}
\]
of $A$.  By the multiplication rules for the components of the radical, together with the assumptions $J_{10}J_{00}=J_{00}J_{01}=\{0\},$
we have $DJ_{00}=J_{00}D=\{0\}.$ Furthermore, $D\cap J_{00}=\{0\}$. Consequently,
$A=D\oplus J_{00}$. Therefore, it is enough to determine
the $T$-ideal of $D$.

Since $J_{10}J_{01}\neq\{0\}$ and  $J_{01}J_{10}=\{0\},$ Lemma~\ref{1} item (7) yields $A_{6}\in\operatorname{var}(D).$
We distinguish two cases.

 \begin{itemize}
    \item \textbf{Case} $C=\{0\}$.
\end{itemize} 

In this case, $J_{11}$ is commutative. Moreover,
\[
[J_{10}J_{01},F]
=[J_{10}J_{01},J_{11}]
=J_{01}J_{10}
=\{0\}.
\]

Set $U=J_{10}$, $V=J_{01}$ and $I=UV$. Since $R$ is commutative, the multiplication rules give
\[
[D,D]^2+D[D,D]D\subseteq I\subseteq Z(D),
\qquad
[D,D]D[D,D]=\{0\},
\]
so $f_4,f_5,f_6\in\operatorname{Id}(D)$.

For $f_1,f_2,f_3$, multilinearity reduces the verification to substitutions from $R$, $U$ and $V$. Evaluations entirely in $R$ vanish. Otherwise, a monomial can be nonzero only with a single $U$-entry at the end, a single $V$-entry at the beginning, or exactly one entry from each, occurring consecutively in the order $UV$. By the commutativity of $R$ and the centrality of $I$, all monomials satisfying
these positional conditions have the same value in each fixed substitution. Their coefficient sum must be zero: replacing the $R$-entries by $e_{11}+e_{33}$ and the
$U$-entries and $V$-entries by $e_{12}$ and $e_{23}$, respectively, gives this same sum multiplying a nonzero matrix unit in $A_6$. Since $f_1,f_2,f_3\in\operatorname{Id}(A_6)$, they vanish on $D$ as well. Hence $\operatorname{Id}(A_6)\subseteq\operatorname{Id}(D)$.

Since $A_{6}\in\operatorname{var}(D)$, the opposite inclusion also
holds. Hence $D\sim_{T}A_{6}.$ It follows that
$A\sim_{T}A_{6}\oplus J_{00}.$

 \begin{itemize}
    \item \textbf{Case} $C\neq \{0\}$.
\end{itemize}

In this case, $R$ is noncommutative and, by the preceding discussion, $R\sim_{T} \mathcal{G}_2.$ Since $R$ is a subalgebra of $D$, we have $\mathcal{G}_2\in\operatorname{var}(D).$ Together with $A_{6}\in\operatorname{var}(D)$, this gives $A_{6}\oplus \mathcal{G}_2\in\operatorname{var}(D),$ and therefore \begin{equation}\label{eq:first-inclusion-A6A2}
\operatorname{Id}(D)
 \subseteq
\operatorname{Id}(A_{6})\cap\operatorname{Id}(\mathcal{G}_2).
\end{equation}

Since $K\notin \textnormal{var}(D)$, by Lemma \ref{commutativityj10j01}, we have $[J_{10}J_{01},J_{11}]=\{0\}$. Thus, setting
$I=J_{10}J_{01}$, we obtain $[I,R]=\{0\}$. Moreover, by \eqref{eq:annihilation-A6},
\[
CI=(CJ_{10})J_{01}=\{0\}
\qquad\text{and}\qquad
IC=J_{10}(J_{01}C)=\{0\}.
\]
Notice also that $P=J_{10}+J_{01}+I$ is a two-sided ideal of $D$. The multiplication rules for the components of the radical, together with $J_{01}J_{10}=\{0\}$, give $P^{2}=I$ and $P^{3}=\{0\}$.

We prove the reverse inclusion in
\eqref{eq:first-inclusion-A6A2}. Let
\[
f=f(x_{1},\ldots,x_{n})
\in \operatorname{Id}(A_{6})\cap\operatorname{Id}(\mathcal{G}_2)
\]
be multilinear. Fix a basis of $D$ obtained as the union of bases of $R$, $J_{10}$ and $J_{01}$. By multilinearity, it is enough to consider evaluations of $f$ on these basis elements.

If every variable is evaluated in $R$, then the resulting value is zero, since $R\sim_{T} \mathcal{G}_2$ and
$f\in\operatorname{Id}(\mathcal{G}_2)$. Furthermore, every evaluation involving at least three entries from $J_{10}\cup J_{01}$
is zero, because $P$ is an ideal of $D$ and $P^{3}=\{0\}$. Evaluations involving exactly two entries both in $J_{10}$ or both in $J_{01}$ also vanish by the multiplication rules.
Thus, it remains to consider evaluations with exactly one entry in $J_{10}\cup J_{01}$, or with one entry in each of $J_{10}$ and $J_{01}$.

Suppose first that exactly one variable, say $x_i$, is evaluated in $J_{10}$. Denote its value by $u$ and the remaining values by $r_{1},\ldots,r_{n-1}\in R$. Since $J_{10}R=\{0\}$, only the monomials in which $x_i$ occurs in the last position can have a nonzero value.
Moreover,
\[
[R,R]J_{10}=CJ_{10}=\{0\},
\] and $RJ_{10}\subseteq J_{10}$. Hence, interchanging adjacent factors from $R$ does not change the value of a product ending in $u$. Consequently, the evaluation of $f$ under consideration is
$\alpha\,r_{1}\cdots r_{n-1}u,$ where $\alpha$ is the sum of the coefficients of the monomials of $f$ in which $x_i$ occurs in the last position.

Consider $a_{6}=e_{11}+e_{33}\in A_{6}$. Evaluating $x_i$ in $e_{12}$ and every other variable in $a_{6}$, we obtain
\[
0=f(a_{6},\ldots,e_{12},\ldots,a_{6})
  =\alpha e_{12},
\]
because $f\in\operatorname{Id}(A_{6})$. Therefore, $\alpha=0$, and the original evaluation vanishes.

Similarly, suppose that exactly one variable $x_i$ is evaluated in $v\in J_{01}$, with all other values in $R$. Since $RJ_{01}=\{0\}$, only the monomials in which $x_i$ occurs in the first position can have a nonzero value. Using
\[
J_{01}[R,R]=J_{01}C=\{0\}
\qquad\text{and}\qquad
J_{01}R\subseteq J_{01},
\] we may reorder the remaining factors without changing their value. Thus, this evaluation is $\beta\,v r_{1}\cdots r_{n-1},$
where $\beta$ is the sum of the coefficients of the monomials in which $x_i$ occurs in the first position. Evaluating $x_i$ in $e_{23}$ and all other variables in $a_{6}$ gives $\beta e_{23}=0$. Hence $\beta=0$, and this evaluation is also zero.

Finally, suppose that $x_i$ is evaluated in
$u\in J_{10}$, $x_j$ is evaluated in $v\in J_{01}$, and all remaining variables are evaluated in
$r_{1},\ldots,r_{n-2}\in R$.
The multiplication rules and $J_{01}J_{10}=\{0\}$ imply that a monomial can have a nonzero value only when $x_i$ occurs immediately before $x_j$.
Indeed, a factor from $R$ immediately after $u$ gives zero, whereas a product in which $v$ precedes $u$ vanishes because $J_{01}R\subseteq J_{01}$ and
$J_{01}J_{10}=\{0\}$.

Therefore, every possibly nonzero monomial has a value of the form
\[
r_{i_{1}}\cdots r_{i_{p}}\,uv\,
r_{j_{1}}\cdots r_{j_{q}},
\qquad p+q=n-2.
\]
Since $uv\in I$ and $[I,R]=\{0\}$, all factors from $R$ can be moved to the left of $uv$. Moreover, $RI\subseteq I$
and
\[
[R,R]I=CI=\{0\},
\]
so these factors can be reordered arbitrarily without changing the value of the product. It follows that the
evaluation of $f$ is $\gamma\,r_{1}\cdots r_{n-2}uv,$ where $\gamma$ is the sum of the coefficients of the monomials in which $x_i$ occurs immediately before $x_j$.

Now evaluate $f$ in $A_{6}$ by setting
$x_i\mapsto e_{12}$, $x_j\mapsto e_{23}$ and $x_k\mapsto a_{6}$ for all $k\notin\{i,j\}.$
Under this evaluation, a monomial has a nonzero value precisely when $x_i$ occurs immediately before $x_j$,
and its value is then $e_{13}$. Since
$f\in\operatorname{Id}(A_{6})$, we obtain
$\gamma e_{13}=0.$ Thus, $\gamma=0$, and the original evaluation vanishes.

We have shown that every evaluation of $f$ in $D$ is zero. Since the field has characteristic zero, this proves
\[
\operatorname{Id}(A_{6})\cap\operatorname{Id}(\mathcal{G}_2)
\subseteq \operatorname{Id}(D).
\]
Together with \eqref{eq:first-inclusion-A6A2}, we obtain
\[
\operatorname{Id}(D)
=
\operatorname{Id}(A_{6})\cap\operatorname{Id}(\mathcal{G}_2)
=
\operatorname{Id}(A_{6}\oplus \mathcal{G}_2).
\]
Therefore, $D\sim_{T}A_{6}\oplus \mathcal{G}_2$. Recalling that $A=D\oplus J_{00}$, it follows that $A\sim_{T}A_{6}\oplus \mathcal{G}_2\oplus J_{00}.$ Since $J_{00}$ is nilpotent, taking $N=J_{00}$ gives the required conclusion in this subcase.
\end{proof}

Recall that $\mathcal J =\{\mathcal G_4,A_7,A_4,A_4^{*},  A_8,A_8^*,\black K\}$ and consider
\[
\mathcal L=\{F,\mathcal G_2,A_1,A_1^*,A_5,A_6\}
\quad \mbox{ and }\quad 
\mathcal B=
\left\{
\bigoplus_{M\in S}M
\mid 
\varnothing\neq S\subseteq\mathcal L
\right\}.
\] 
\begin{theorem}\label{ClassificationLinear}
Let $A$ be an algebra over a field $F$ of characteristic zero. Then, the following conditions are equivalent:
\begin{enumerate}
   \item There exists a constant $\alpha>0$ such that
$c_n^z(A)\leq \alpha n$ for every $n\geq 1$.

\item $A$ is PI-equivalent to either $N$ or $B\oplus N$, where $B\in\mathcal B$ and $N$ is a nilpotent algebra.
    \item $Q\notin \textnormal{var}(A)$, for all $Q\in \mathcal{J}$.
\end{enumerate}
\end{theorem}

    \begin{proof}
By Lemmas~\ref{lem:UT3-central}--\ref{G2k} and \ref{k}, the central codimensions of every
$Q\in\mathcal{J}$ are not bounded by a linear function. Since $c_n^z(Q)\leq c_n^z(A)$ whenever $Q\in\operatorname{var}(A)$, this
proves $(1)\Rightarrow(3)$.

Moreover, $c_n^z(F)=0$, and Lemmas \ref{A1A1*}--\ref{G2k} show that every algebra in $\mathcal{L}$ has central codimensions bounded
by a linear function. By Remark~\ref{directsum}, the same holds for their finite direct sums. Thus, we obtain $(2)\Rightarrow(1)$.

It remains to prove that $(3)\Rightarrow(2)$.
Assume that $Q\notin\operatorname{var}(A)$, for all
$Q\in\mathcal{J}$. Since
$\mathcal{G}_4\in\operatorname{var}(\mathcal{G})$ and
$A_4\in\operatorname{var}(UT_2)$, we also have
$\mathcal{G},UT_2\notin\operatorname{var}(A)$.
Hence, by Proposition~\ref{teokemer}, after replacing $A$ by a PI-equivalent algebra, we may assume that
$$
A=B_1\oplus\cdots\oplus B_m,
$$ where each $B_i$ is finite-dimensional and is either nilpotent or satisfies $B_i=F\dotplus J(B_i)$. The case in which all the summands are nilpotent is immediate.
Thus, let $B_i$ be a non-nilpotent summand. As before, write
\[
J(B_i)=J_{10}\dotplus J_{01}\dotplus J_{00}\dotplus J_{11}.
\]
Since $B_i\in\operatorname{var}(A)$ and
$A_4,A_4^*\notin\operatorname{var}(A)$, items~(4) and~(5)
of Lemma~\ref{1} imply that $J_{10}J_{00}=J_{00}J_{01}=\{0\}.$ Set
\[
D_i=F\dotplus J_{10}\dotplus J_{01}\dotplus J_{11}\dotplus J_{01}J_{10}.
\]
By the multiplication rules for the components of the Jacobson radical, $D_i$ is a subalgebra of $B_i$. Moreover,
$D_iJ_{00}=J_{00}D_i=\{0\}.$ Although $D_i\cap J_{00}=J_{01}J_{10}$ may be nonzero, every mixed product involving elements of $D_i$ and $J_{00}$ vanishes.
Therefore, by multilinearity, a multilinear polynomial is an identity of $B_i=D_i+J_{00}$ if and only if it is an identity of both $D_i$ and $J_{00}$. Consequently, since the field has characteristic zero,
\[
\operatorname{Id}(B_i)
=
\operatorname{Id}(D_i)\cap\operatorname{Id}(J_{00})
=
\operatorname{Id}(D_i\oplus J_{00}),
\]
and hence $B_i\sim_T D_i\oplus N_i$, where
$N_i=J_{00}$ is nilpotent.

We now reduce the study of $D_i$ to two simpler cases. Define
\[
I_5=J_{01}J_{10}
\qquad\text{and}\qquad
I_6=J_{10}J_{01}.
\]
The multiplication rules and the preceding vanishing relations show that both $I_5$ and $I_6$ are two-sided ideals of $D_i$.
Moreover, $I_5\cap I_6\subseteq J_{00}\cap J_{11}=\{0\}.$ Consider the canonical homomorphism
\[
\Phi_i:D_i\rightarrow D_i/I_5\oplus D_i/I_6,
\qquad
\Phi_i(d)=(d+I_5,d+I_6).
\]
Since $\ker\Phi_i=I_5\cap I_6=\{0\}$, the map $\Phi_i$ is
injective, and $D_i$ is isomorphic to a subalgebra of $D_i/I_5\oplus D_i/I_6$.
Conversely, both quotients are homomorphic images of $D_i$. Since varieties are closed under finite direct sums, we conclude
that $D_i\sim_T D_i/I_5\oplus D_i/I_6.$

It is therefore enough to study the two quotient algebras separately. Both belong to $\operatorname{var}(A)$, so all the exclusions in condition~(3) hold for them.
Notice that the $00$-component of $J(D_i)$ is precisely $I_5$. In what follows, we use the same notation $J_{rs}$ for the radical components of the quotient being considered.

In $D_i/I_5$, we have
$J_{00}=J_{01}J_{10}=\{0\},$
and so
$D_i/I_5=F\dotplus J_{10}\dotplus J_{01}\dotplus J_{11}$. Set
\[
P_i=D_i/I_5
=
F\dotplus J_{10}\dotplus J_{01}\dotplus J_{11},
\]
where the notation $J_{rs}$ now refers to the images of the corresponding radical components in the quotient.

If $J_{10}J_{01}\neq\{0\}$, then Lemma~\ref{lemma-A6} gives either
$P_i\sim_T A_6\oplus N$ or  $P_i\sim_T A_6\oplus \mathcal{G}_2\oplus N,$ for some nilpotent algebra $N$.

If $J_{10}J_{01}=\{0\}$, then Lemma~\ref{1} yields $P_i\in\operatorname{var}(UT_2).$
Recalling that $A_4=A_{3,1}$, $A_4^*=A_{3,1}^*$,
$A_7\sim_T N_4$ and $\mathcal{G}_2\sim_T N_3$, the exclusions in condition~(3), together with Theorem~\ref{subvarieties}, show that $P_i$ is $T$-equivalent to a nonempty direct sum of some of the algebras in $\{F,\mathcal{G}_2,A_1, A_1^*\}$, possibly together with a nilpotent direct summand.

Similarly, in $D_i/I_6$, we have
\[
J_{10}J_{00}=J_{00}J_{01}=J_{10}J_{01}=\{0\},
\]
and its $00$-component is $J_{00}=J_{01}J_{10}.$  If $J_{01}J_{10}\neq\{0\}$, then Lemma~\ref{lemma-A5} implies that either $D_i/I_6\sim_T A_5\oplus N$ or $D_i/I_6\sim_T A_5\oplus\mathcal G_2\oplus N,$
for some nilpotent algebra $N$.

Assume that $J_{01}J_{10}=\{0\}$. Then $J_{00}=\{0\}$ and $D_i/I_6=R_i$ where $R_i=F\dotplus J_{10}\dotplus J_{01}\dotplus J_{11}.$ Since also $J_{10}J_{01}=\{0\}$, Lemma~\ref{1} yields $R_i\in\operatorname{var}(UT_2).$ By the same argument used for $P_i$, Theorem~\ref{subvarieties} implies that $R_i$ is $T$-equivalent to $B\oplus N$, where $N$ is nilpotent and $B$ is a finite direct sum of algebras from $\{F,\mathcal G_2,A_1,A_1^*\}$.

Finally, recall that $A=B_1\oplus\cdots\oplus B_m,$ where each $B_i$ is either nilpotent or satisfies
\[
B_i\sim_T D_i/I_5\oplus D_i/I_6\oplus N_i.
\]

Taking the direct sum of all these components, collecting all nilpotent summands, and omitting repeated summands, we obtain $A\sim_T B\oplus N$, where $B$ is a finite direct sum of algebras
from $\mathcal L$, completing the proof.
\end{proof}

\begin{remark}\label{rem:central-incomparability}
The ordinary $T$-ideals of the algebras in
$\{\mathcal G_4,A_7,A_4,A_4^*,K\}$
are pairwise incomparable.
Indeed, let $St_5=\sum_{\sigma\in S_5}(\textnormal{sgn}\, \sigma)\,  x_{\sigma(1)}\cdots x_{\sigma(5)}$ denote the standard polynomial of degree five. A direct calculation gives the following identities:
\[
\begin{array}{c|l}
\text{Algebra} & \text{Polynomial identities}\\ \hline
\mathcal G_4 & [x_1,x_2,x_3]\\
A_7 & [x_1,x_2,x_3,x_4],\,\, [x_1,x_2][x_3,x_4]\\
A_4 & [x_1,x_2]x_3x_4\\
A_4^* & x_3x_4[x_1,x_2]\\
K & g=[x_1[x_1,x_2]x_1,x_3],\,\, St_5.
\end{array}
\]
Every other member of $\mathcal J-\{A_8,A_8^*\}$ fails at least one polynomial
in the corresponding row.  Thus, for distinct $Q,R\in\mathcal J-\{A_8,A_8^*\}$,
$\operatorname{Id}(Q)\not\subseteq\operatorname{Id}(R).$ Moreover, by monotonicity of central codimensions, $A_8,A_8^*\notin \textnormal{var}(Q)$, for all $Q\in \{\mathcal G_4, A_7, A_4, A_4^*, K\}.$
\end{remark}

\begin{corollary}
{ The $z$-minimal varieties of  degree 2}
over a field of characteristic zero are precisely the varieties
generated by $\mathcal G_4$, $A_7$, $A_4$, $A_4^*$, $K$.
\end{corollary}

\begin{proof}
Suppose that $\mathcal V$ is { $z$-minimal of degree 2.} Since $\mathcal V$ has quadratic central codimension growth and $A_8$ and $A_8^*$ have cubic central codimension growth, it follows that $A_8,A_8^*\notin \mathcal V$. Moreover, by Theorem~\ref{ClassificationLinear},
some $Q\in\mathcal J-\{A_8,A_8^*\}$ belongs to $\mathcal V$.
Since $\operatorname{var}(Q)\subseteq\mathcal V$ and both
varieties have quadratic central codimension growth,
$z$-minimality gives $\mathcal V=\operatorname{var}(Q)$.

Conversely, fix $Q\in\mathcal J-\{A_8,A_8^*\}$ and let
$\mathcal U\subsetneq\operatorname{var}(Q)$.
Then $Q\notin\mathcal U$.
By Remark~\ref{rem:central-incomparability}, no other member
of $\mathcal J$ belongs to $\operatorname{var}(Q)$, and hence
$\mathcal J\cap\mathcal U=\varnothing$.
Theorem~\ref{ClassificationLinear} now gives
$c_n^z(\mathcal U)\leq \alpha n$, proving the required $z$-minimality.
\end{proof}

\section{Applications in the graded setting}
\label{sec:graded-bounded}

In this section, we apply the classification of ordinary algebras with bounded central codimensions obtained in Section~4 to the graded setting. More precisely, for a finite group $G$, we characterize the $G$-graded
algebras whose graded central codimension sequences are bounded by a constant. The classification is given up to $T_G$-equivalence and permits an arbitrary commutative graded summand.

Throughout this section, $G$ is a finite group with identity element
$1_G$. A $G$-graded algebra is an algebra endowed with a decomposition
\[
A=\bigoplus_{g\in G}A_g
\]
into $F$-subspaces such that $A_gA_h\subseteq A_{gh}$ for all $g,h\in G$.
The subspace $A_g$ is called the homogeneous component of degree $g$ and for $0\neq a\in A_g$, we write $\deg(a)=g$.
A subalgebra $B$ of $A$ is graded if
\[
B=\bigoplus_{g\in G}(B\cap A_g).
\]

If $A$ is unitary, then $1_A\in A_{1_G}$. Let $F\langle X,G\rangle$ be the free associative
$G$-graded algebra without constant term on the disjoint sets
$X_g=\{x_{1,g},x_{2,g},\ldots\}$, where $\deg x_{i,g}=g$.
We denote by $\operatorname{Id}^G(A)$ the $T_G$-ideal of graded
polynomial identities of $A$, and by $\operatorname{Id}^{G,z}(A)$ the
$T_G$-space of polynomials whose graded evaluations belong to $Z(A)$.
Here $Z(A)$ is the ordinary center of $A$, and graded identities are included among graded central polynomials. Moreover, a $T_G$-ideal, respectively a $T_G$-space, is an ideal, respectively an $F$-subspace, invariant under all graded endomorphisms
of the free algebra.

The graded variety $\operatorname{var}^G(A)$ generated by $A$ is the class of all $G$-graded algebras $B$ such that
$
\operatorname{Id}^G(A)\subseteq\operatorname{Id}^G(B).$ Two $G$-graded algebras $A$ and $B$ are $T_G$-equivalent, written
$A\sim_{T_G}B$, if $\operatorname{Id}^G(A)=\operatorname{Id}^G(B)$.

For an ordinary algebra $R$, we use the same notation for $R$ when it is endowed with the trivial $G$-grading, namely,
$R_{1_G}=R$ and  $R_g=\{0\}$ for all $g\neq 1_G.$

For $n\geq 1$, set
\[
P_{n}^G
 =\operatorname{span}_F
 \{x_{\sigma(1),g_{\sigma(1)}}\cdots
   x_{\sigma(n),g_{\sigma(n)}}\mid \sigma\in S_n, g_i\in G\}.
\]
The corresponding codimensions are
\[
c_n^G(A)
 =\dim_F\frac{P_n^G}{P_n^G\cap\operatorname{Id}^G(A)},
\qquad
c_n^{G,z}(A)
 =\dim_F\frac{P_n^G}{P_n^G\cap\operatorname{Id}^{G,z}(A)}.
\]

A graded polynomial $f$ is central on $A$ if and only if
$[f,x_{j,h}]\in\operatorname{Id}^G(A)$ for every $h\in G$, where the
variables $x_{j,h}$ do not occur in $f$. Thus graded central
polynomials are determined by graded identities. Consequently,
$T_G$-equivalent algebras have the same graded central codimensions,
and $B\in\operatorname{var}^G(A)$ implies $
c_n^{G,z}(B)\leq c_n^{G,z}(A).$
Moreover, since
$\operatorname{Id}^{G,z}(A\oplus B)
 =\operatorname{Id}^{G,z}(A)\cap\operatorname{Id}^{G,z}(B)$,
\begin{equation*}\label{eq:graded-central-sums}
\max\{c_n^{G,z}(A),c_n^{G,z}(B)\}
\leq c_n^{G,z}(A\oplus B)
\leq c_n^{G,z}(A)+c_n^{G,z}(B).
\end{equation*}
If $B^q=\{0\}$, then
$c_n^{G,z}(A\oplus B)=c_n^{G,z}(A)$ for all $n\geq q$.

For every $G$-graded algebra $A$, regarded also as an ordinary algebra,
we have $c_n^z(A)\leq c_n^{G,z}(A)$, for all $n\geq 1.$ Indeed, consider the linear map
\[
\iota_n:P_n\longrightarrow P_n^G,
\qquad
\iota_n(f)=
f\left(\sum_{g\in G}x_{1,g},\ldots,
       \sum_{g\in G}x_{n,g}\right).
\]
Every element of $A$ is a sum of homogeneous elements, so $\iota_n(f)\in\operatorname{Id}^{G,z}(A)$ if and only if $f\in\operatorname{Id}^{z}(A).$
Hence $\iota_n$ induces an injective linear map
\[
\frac{P_n}{P_n\cap\operatorname{Id}^{z}(A)}
\hookrightarrow
\frac{P_n^G}{P_n^G\cap\operatorname{Id}^{G,z}(A)},
\]
which proves the inequality.

Let $L$ be a field extension of $F$. For any $G$-graded $L$-algebra $E$, let $E_F$ denote the same algebra regarded as an $F$-algebra by restriction of scalars. Since this operation does not change the multiplication,
the homogeneous components, or the center, we have
\[
\operatorname{Id}_F^G(E_F)
=
\operatorname{Id}_L^G(E)\cap F\langle X,G\rangle
\quad \mbox{ and }\quad  \operatorname{Id}_F^{G,z}(E_F)
=
\operatorname{Id}_L^{G,z}(E)\cap F\langle X,G\rangle.
\]
Recall also that, for every $G$-graded $F$-algebra $R$,
\[
\operatorname{Id}_F^G(R\otimes_F L)
=
\operatorname{Id}_F^G(R)
\qquad\text{and}\qquad
\operatorname{Id}_F^{G,z}(R\otimes_F L)
=
\operatorname{Id}_F^{G,z}(R).
\]
Indeed, a graded polynomial is central if and only if its commutator with a new homogeneous variable of each degree is a graded identity.

Consequently,
\[c_{n,L}^{G}(R\otimes_F L)
=
c_{n,F}^{G}(R)\quad  \mbox{ and }\quad 
c_{n,L}^{G,z}(R\otimes_F L)
=
c_{n,F}^{G,z}(R)
\qquad\text{for every }n\geq 1,
\]
where the subscripts indicate the fields over which the
dimensions are computed.

\black

Unlike the ordinary case, bounded graded central codimensions do not imply polynomial growth of the graded codimensions.

\begin{example} Fix $g\in G-\{1_G\}$ of order $m> 1$ and consider $H=\langle g\rangle$.
\begin{enumerate} 
    \item Endow the commutative group algebra $FH$ with its natural $G$-grading given by $(FH)_h=Fh$ for $h\in H$ and zero components outside $H$. Then $c_n^G(FH)=m^n$ for all $n\geq 1$, and so $FH$ has exponential $G$-graded codimension growth.
 
  \item   The commutative unitary
algebra $C_2=F(e_{11}+e_{22})+Fe_{12}
\subseteq UT_2$, graded by $\deg (e_{11}+e_{22})=1_G$ and $\deg(e_{12})=g$
satisfies $c_n^G(C_2)=n+1.$ \end{enumerate}

Both algebras are commutative and therefore satisfy $c_n^{G,z}(FH)=c_n^{G,z}(C_2)=0$
for all $n\geq1.$
\end{example}

More generally, every commutative $G$-graded algebra has zero graded central codimensions, since every graded polynomial takes its values
in the center. Thus bounded graded central codimensions impose no restriction on the graded identities of a commutative direct summand. The classification below reflects this fact by allowing an arbitrary commutative graded component.

 For $g\in G$, let
\[
A_1^g=Fe_{11}+Fe_{12}
\quad \mbox{ and }\quad 
(A_1^*)^g=Fe_{22}+Fe_{12},
\]
with $\deg e_{11}=\deg e_{22}=1_G$ and $\deg e_{12}=g$.
For $g=1_G$, these are the trivially graded versions of $A_1$ and
$A_1^*$ introduced in Section~3.

\begin{lemma} \cite[Theorem 5.4]{KoshlukovLaMattina2015} \label{lem:graded-triangular-central}
For every $g\in G$ and $n\geq1$,
\[
c_n^{G,z}(A_1^g)=c_n^G(A_1^g)
=c_n^{G,z}((A_1^*)^g)=c_n^G((A_1^*)^g)
=\begin{cases}
n,&g=1_G,\\
n+1,&g\neq1_G.
\end{cases}
\]
\end{lemma}

The following lemma provides a structural decomposition of the algebras under consideration.

\begin{lemma}\label{lem:finite-dimensional-bounded-central}
Let $A=\bigoplus_{g\in G}A_g$ be a finite-dimensional $G$-graded algebra
over a field $F$ of characteristic zero. Suppose that $c_n^{G,z}(A)\leq \alpha$  
for all $n\geq 1$ and for a constant $\alpha$.
Then $A$ admits a decomposition
\[
A=U\oplus N
\]
as a direct sum of graded algebras, where $N$ is nilpotent and $U$ is either zero or a unitary algebra.
\end{lemma}

\begin{proof}
By the graded Wedderburn--Malcev theorem, write $A=B\dotplus J(A),$ where $B$ is a graded semisimple subalgebra. If $B=\{0\}$, then
$A=J(A)$ is nilpotent, and the result follows by taking
$U=\{0\}$ and $N=A$.

Suppose that $B\neq\{0\}$ and set $e=1_B$. Since
$e\in B_{1_G}$, the decomposition
\[
J(A)=J_{00}\dotplus J_{01}\dotplus J_{10}\dotplus J_{11},
\qquad
J_{ij}=\{a\in J(A)\mid ea=ia,\ ae=ja\},
\]
is a decomposition into graded subspaces. If $J_{10}\neq\{0\}$, choose a nonzero homogeneous element
$a\in J_{10}\cap A_g$. Then
\[
e^2=e,\qquad ea=a,\qquad ae=a^2=0.
\]
Thus $Fe+Fa$ is a graded subalgebra of $A$ isomorphic to $A_1^g$.
By Lemma~\ref{lem:graded-triangular-central},
$c_n^{G,z}(A)\geq c_n^{G,z}(A_1^g)\geq n$
for all $n\geq 1,$ contrary to the hypothesis. Similarly, a nonzero homogeneous element
of $J_{01}$ would yield a graded subalgebra isomorphic to $(A_1^*)^g$.
Hence $J_{10}=J_{01}=\{0\}.$

Set $U=B\dotplus J_{11}$ and $N=J_{00}.$
The multiplication rules give $UN=NU=\{0\}$, so
$A=U\oplus N$ is a direct sum of graded algebras.
Moreover, $N$ is nilpotent and $e$ is the identity element of $U$.
\end{proof}

In the decomposition of Lemma~\ref{lem:finite-dimensional-bounded-central}, the nilpotent summand $N$ does not affect the central codimensions in sufficiently large degrees. 
Therefore, up to a nilpotent graded summand, it suffices to study
finite-dimensional unitary $G$-graded algebras $U$ satisfying
$c_n^{G,z}(U)\leq\alpha$ for all $n\geq1$ and some constant $\alpha>0$.

\begin{lemma}\label{lem:nonneutral-components-central}
Let $A=\bigoplus_{g\in G}A_g$ be a unitary $G$-graded algebra over a field $F$ of characteristic zero. If there exists a constant $\alpha >0$ such that $c_n^{G,z}(A)\leq \alpha$ for all $n\geq 1$, then $A_gA_{1_G}\subseteq Z(A)$ for every $g\in G-\{1_G\}.$ In particular, $A_g\subseteq Z(A)$ for every $g\in G-\{1_G\}.$
\end{lemma}

\begin{proof}
Fix $g\in G-\{1_G\}$ and suppose, by contradiction, that there exist
$a\in A_g$ and $b\in A_{1_G}$ such that $ab\notin Z(A)$.
For each $n\geq 2$, consider the multilinear graded monomials
\[
m_i=x_{i,g}x_{1,1_G}\cdots
\widehat{x_{i,1_G}}\cdots x_{n,1_G},
\qquad 1\leq i\leq n,
\]
where $\widehat{x}$ indicates that the corresponding variable is omitted.

We claim that the classes of $m_1,\ldots,m_n$ are linearly independent modulo ${P_n^G\cap\operatorname{Id}^{G,z}(A)}.$ Indeed, suppose that
\[
f=\sum_{k=1}^{n}\lambda_km_k
\in\operatorname{Id}^{G,z}(A),
\qquad \lambda_k\in F.
\]
Fix $i\in\{1,\ldots,n\}$ and choose
$j\in\{1,\ldots,n\}\setminus\{i\}$.
Since $1_A\in A_{1_G}$, consider the graded evaluation
\[
\begin{aligned}
x_{i,g}\mapsto a,
\quad 
x_{t,g}\mapsto 0,\quad 
x_{j,1_G}\mapsto b,\quad  x_{l,1_G}\mapsto 1_A,
\end{aligned}
\] for all $t\neq i$ and $l\neq j.$ Under this evaluation, $m_i$ takes the value $ab$, whereas every
$m_k$ with $k\neq i$ vanishes. Hence $f$ takes the value
$\lambda_iab$. Since $f$ is a graded central polynomial of $A$, we obtain $\lambda_iab\in Z(A).$
As $ab\notin Z(A)$, it follows that $\lambda_i=0$.
Since $i$ was arbitrary, all coefficients vanish, proving the claim.

Consequently, $c_n^{G,z}(A)\geq n$ for all  $n\geq 2$, a contradiction. Therefore,
$A_gA_{1_G}\subseteq Z(A)$ for every $g\neq 1_G$.

Finally, since $1_A\in A_{1_G}$, we have
\[
A_g=A_g1_A\subseteq A_gA_{1_G}\subseteq Z(A),
\]
which proves the last assertion.
\end{proof}

\begin{lemma}\label{lem:neutral-component-and-center}
Let $A=\bigoplus_{g\in G}A_g$ be a unitary $G$-graded algebra over a field $F$ of characteristic zero. Suppose that
$A_g\subseteq Z(A)$ for every $g\in G-\{1_G\}.$ Then
\[
A\sim_{T_G}A_{1_G}\oplus Z(A),
\]
where $A_{1_G}$ is endowed with the trivial $G$-grading and $Z(A)$
with the grading induced by $A$.
\end{lemma}

\begin{proof}
Since every nonneutral homogeneous component of $A$ is central,
\[
Z(A)=\bigl(Z(A)\cap A_{1_G}\bigr)
\dotplus (\sum_{g\neq1_G}A_g).
\]
Thus $Z(A)$ is a graded subalgebra of $A$. As $A_{1_G}$ is also
a graded subalgebra, we have
\[
\operatorname{Id}^G(A)
\subseteq
\operatorname{Id}^G(A_{1_G})
\cap\operatorname{Id}^G(Z(A))
=
\operatorname{Id}^G(A_{1_G}\oplus Z(A)).
\]

To prove the reverse inclusion, first observe that, for
$c\in A_g$ with $g\neq1_G$ and $a,b\in A_{1_G}$, both $c$ and
$ca\in A_g$ are central. Hence
\[0=[ca,b]=c[a,b]+[c,b]a=c[a,b].
\]
Consequently,
\[
A_g[A_{1_G},A_{1_G}]=0
\qquad\text{for every }g\neq1_G.
\]

By multihomogeneity and multilinearization, it suffices to consider a multilinear
graded polynomial
\[
f=f(y_1,\ldots,y_r,z_1,\ldots,z_s)
\in
\operatorname{Id}^G(A_{1_G})
\cap\operatorname{Id}^G(Z(A)),
\]
where each $y_i$ has degree $1_G$ and each $z_j$ has degree
$g_j\neq1_G$.

If $s=0$, every graded evaluation of $f$ in $A$ takes place in
$A_{1_G}$, and therefore vanishes. Suppose that $s\geq1$, and
consider an arbitrary graded evaluation
$y_i\mapsto b_i\in A_{1_G}$ and  $z_j\mapsto c_j\in A_{g_j}.$
All the $c_j$ are central, so they can be moved past every factor.
Moreover, since $c_1[A_{1_G},A_{1_G}]=0$, the neutral factors
can be permuted in any product containing $c_1$. Checking each
monomial of $f$, we therefore obtain
\[
f(b_1,\ldots,b_r,c_1,\ldots,c_s)
=
b_1\cdots b_r\,
f(1_A,\ldots,1_A,c_1,\ldots,c_s),
\]
where the empty product is understood to be $1_A$.

Since $1_A\in Z(A)\cap A_{1_G}$ and
$c_j\in Z(A)\cap A_{g_j}$ for every $j$, the evaluation on the
right vanishes since $f\in\operatorname{Id}^G(Z(A))$.
Thus $f\in\operatorname{Id}^G(A)$, proving
\[
\operatorname{Id}^G(A)
=
\operatorname{Id}^G(A_{1_G})
\cap\operatorname{Id}^G(Z(A)).
\]
The desired $T_G$-equivalence follows.
\end{proof}

We are now ready to prove the main result of this section.

\begin{theorem}\label{thm:graded-bounded-classification}
Let $A$ be a finite-dimensional $G$-graded algebra over a field $F$ of characteristic
zero. Then the following conditions are equivalent:
\begin{enumerate}
 \item There exists a constant $\alpha>0$ such that $c_n^{G,z}(A)\leq\alpha$ for all $n\geq1.$
    \item There exist a finite-dimensional commutative $G$-graded algebra $C$, which is either zero or unital, and a finite-dimensional nilpotent $G$-graded algebra $N$ such that either
    \[
    A\sim_{T_G} C \oplus  N \quad \mbox{ or }\quad  A\sim_{T_G} \mathcal{G}_2\oplus C \oplus  N
    \]
    where $\mathcal G_2$ is endowed with the trivial $G$-grading.
\end{enumerate}
\end{theorem}

\begin{proof}
Suppose first that condition~(2) holds, and choose $s\geq1$ such
that $N^s=0$. Since $C$ is commutative, every graded polynomial
is central on $C$. Moreover, every multilinear graded polynomial
of degree at least $s$ is an identity of $N$. Therefore, for all
$n\geq s$,
\[
c_n^{G,z}(A)=
\begin{cases}
0, & \text{if } A\sim_{T_G}C\oplus N,\\
1, & \text{if } A\sim_{T_G}\mathcal G_2\oplus C\oplus N,
\end{cases}
\]
where the second equality follows from Lemma~\ref{G2k} and the trivial
grading on $\mathcal G_2$. Since the remaining codimensions are
finite, condition~(1) follows.

Conversely, suppose that $c_n^{G,z}(A)\leq\alpha$ for all $n\geq1$.
By Lemma~\ref{lem:finite-dimensional-bounded-central}, we may write $A=U\oplus N,$ where $N$ is a finite-dimensional nilpotent graded algebra and
$U$ is either zero or a finite-dimensional unitary graded algebra.
If $U=\{0\}$, the conclusion follows by taking $C=\{0\}$.

Suppose that $U\neq0$. Since $c_n^{G,z}(U)\leq c_n^{G,z}(A)\leq\alpha,$ for all $n\geq1,$ Lemma~\ref{lem:nonneutral-components-central} gives $U_g\subseteq Z(U)$ for every $g\neq1_G$.
Thus, by Lemma~\ref{lem:neutral-component-and-center}, $U\sim_{T_G}U_{1_G}\oplus Z(U),$
where $U_{1_G}$ has the trivial grading and $Z(U)$ has the grading induced by $U$. Set $C=Z(U)$. Then $C$ is a finite-dimensional commutative unitary  graded algebra and
\[
A\sim_{T_G}U_{1_G}\oplus C\oplus N.
\]

The algebra $U_{1_G}$ is unitary and satisfies
\[
c_n^z(U_{1_G})
=
c_n^{G,z}(U_{1_G})
\leq c_n^{G,z}(U)
\leq\alpha.
\]

{By Theorem~\ref{boundedbyconstant}, since
$U_{1_G}$ is nonzero and unital, we have
$U_{1_G}\sim_T D\oplus M$,
where $D\in\{F,\mathcal G_2\}$ and
$M$ is nilpotent.
By Remark~\ref{rem:unital-nilpotent},
we obtain $U_{1_G}\sim_T D$.
Since both algebras are trivially graded,
it follows that $U_{1_G}\sim_{T_G}D$.}

Hence $A\sim_{T_G}D\oplus C\oplus N.$ If $D=\mathcal G_2$, this is the required conclusion. If $D=F$, observe that $F1_U$ is a graded subalgebra of $C=Z(U)$.
Consequently,
$\operatorname{Id}^G(C)\subseteq\operatorname{Id}^G(F),$ and therefore $F\oplus C\sim_{T_G}C$. This proves condition~(2). \end{proof}

Combining Theorem~\ref{thm:graded-bounded-classification} with the graded representability results of \cite{AljadeffKanelBelov2010},
we now remove the finite-dimensional assumption. 

\begin{corollary}
  Let $G$ be a finite group and let $A$ be a $G$-graded algebra over a field $F$ of characteristic zero. Then $c_n^{G,z}(A)\leq\alpha$ for all $n\geq1$ and some constant $\alpha>0$ if and only if $A$ is $T_G$-equivalent to either
  $$C\oplus N \quad \mbox{ or }\quad  \mathcal G_2\oplus C\oplus N,$$ for some commutative $G$-graded algebra $C$, which is either zero
or unital, and a nilpotent $G$-graded algebra $N$, where $\mathcal G_2$ is endowed with the trivial $G$-grading.
\end{corollary}

\begin{proof}
Suppose first that $c_n^{G,z}(A)\leq\alpha$ for all $n\geq1$.
The case $A=\{0\}$ is immediate, so assume that $A\neq \{0\}$.
Since
\[
c_n^z(A)\leq c_n^{G,z}(A)\leq\alpha,
\]
Proposition~\ref{teokemer} provides a finite-dimensional
ordinary algebra $D_0$ such that
$\operatorname{Id}(A)=\operatorname{Id}(D_0).$ Set $d=\dim_F D_0$. Every multilinear graded polynomial
alternating in $d+1$ variables of the same homogeneous degree
is therefore an identity of $A$.

By the argument in
\cite[Section~7.1]{AljadeffKanelBelov2010}, the relatively free
$G$-graded algebra $R$ on $d$ generators of each degree in
$\operatorname{var}^G(A)$ has the same graded identities as $A$.
Since $G$ is finite, $R$ is affine. Moreover, $R$ satisfies
the ordinary identities of $A$, and hence is a PI-algebra.
Thus \cite[Theorem~1.1]{AljadeffKanelBelov2010} yields a field
extension $L/F$ and a finite-dimensional $G$-graded $L$-algebra
$D$ such that $\operatorname{Id}_F^G(D)=\operatorname{Id}^G(A)$.

Since graded identities determine graded central polynomials, $\operatorname{Id}_F^{G,z}(D)
=
\operatorname{Id}^{G,z}(A).$ Furthermore, the $L$-linear span of the graded central
polynomials of $D$ with coefficients in $F$ consists of
graded central polynomials with coefficients in $L$.
Consequently,
\[
c_{n,L}^{G,z}(D)
\leq c_{n,F}^{G,z}(D)
=c_n^{G,z}(A)
\leq\alpha.
\]

Applying Theorem~\ref{thm:graded-bounded-classification} over $L$, we obtain a commutative $G$-graded $L$-algebra $C'$, which is either zero or unital, a nilpotent $G$-graded $L$-algebra $N'$, and
$B'\in\{\{0\},\mathcal G_2(L)\}$ such that
\[
\operatorname{Id}_L^G(D)
=
\operatorname{Id}_L^G(B'\oplus C'\oplus N').
\]

Set $C=C'_F$ and $N=N'_F$, where the subscript $F$ denotes
restriction of scalars. Then $C$ is a commutative $G$-graded
$F$-algebra and $N$ is a nilpotent $G$-graded $F$-algebra.
Choose $B=0$ if $B'=0$ and $B=\mathcal G_2(F)$ otherwise,
endowed with the trivial grading. Since
$B'\cong B\otimes_F L$, the invariance of graded identities
under scalar extension gives $\operatorname{Id}_F^G(B')=\operatorname{Id}^G(B).$ Intersecting the equality of $L$-identities above with
$F\langle X,G\rangle$, we therefore obtain
\[
\begin{aligned}
\operatorname{Id}^G(A)=\operatorname{Id}_F^G(D)=\operatorname{Id}_F^G(B')
  \cap\operatorname{Id}^G(C)\cap\operatorname{Id}^G(N)=\operatorname{Id}^G(B\oplus C\oplus N).
\end{aligned}
\]
Thus $A\sim_{T_G}B\oplus C\oplus N$, as required.

Conversely, suppose that
$A\sim_{T_G}B\oplus C\oplus N$, where
$B\in\{\{0\},\mathcal G_2\}$ and $C,N$ have the stated properties.
Choose $s\geq1$ such that $N^s=\{0\}$.
Since $C$ is commutative and every multilinear graded
polynomial of degree at least $s$ is an identity of $N$,
we obtain
\[
c_n^{G,z}(A)=c_n^{G,z}(B)=
\begin{cases}
0,& B=\{0\},\\
1,& B=\mathcal G_2,
\end{cases}
\qquad\text{for all }n\geq s.
\]
The second case follows from Lemma~\ref{G2k} and the trivial
grading on $\mathcal G_2$. Since $G$ is finite, the remaining
codimensions are finite, and hence the sequence
$c_n^{G,z}(A)$ is bounded by a constant.
\end{proof}

\section*{AI Declaration}

During the preparation of this manuscript, the authors used
OpenAI's ChatGPT (GPT-6) to assist with language editing, to help identify relevant references and to check selected mathematical arguments. The authors independently verified the mathematical content and take full responsibility for the manuscript.


\begin{thebibliography}{999}

\bibitem{AljadeffKanelBelov2010}
E. Aljadeff and A. Kanel-Belov.
{\it Representability and Specht problem for $G$-graded algebras.}
Adv. Math. 225 (2010) 2391--2428.


\bibitem{Berele}
A. Berele.
{\it Properties of hook Schur functions with applications to p.\,i. algebras}.
{Adv. in Appl. Math.} {41} (2008), 52--75.

\bibitem{BR}
A. Berele and A. Regev.
{\it Asymptotic behaviour of codimensions of p.\,i. algebras satisfying Capelli identities}.
{Trans. Amer. Math. Soc.} {360} (2008), 5155--5172.

\bibitem{AnaThaisWillerJuan} W. D. S. Costa, J. P. Cruz, T. S. do Nascimento and A. C. Vieira. \textit{ Central cocharacters of the subvarieties of varieties of superalgebras with almost polynomial growth}.
arXiv:2511.10511 (2025).

\bibitem{CMQ}
W. Q. Cota, L. H. de Souza Matos and P. Quintino.
Varieties of associative algebras with quadratic codimension growth.
arXiv:2608.22966v2 (2026).

\bibitem{NascimentoCota} W. Q. Cota and T. S. do Nascimento. \newblock{\em On the multiplicities of the central cocharacter of algebras with polynomial identities}. J. Algebra 697 (2026) 11--37.

\bibitem{Drensky} V. Drensky. Free Algebras and PI-Algebras: Graduate Course in Algebra.
Springer-Verlag, Singapore, 2000.


\bibitem{GLa} A. Giambruno and D. La Mattina.  
		\newblock{\em PI-algebras with slow codimension growth.} \newblock J. Algebra 284 (2005) 371--391.



 \bibitem{GiamZai} A. Giambruno and M. Zaicev. \newblock{\em A characterization of algebras with polynomial growth of
the codimensions}. Proc. Amer. Math. Soc. 129 (2001) 59--67.


\bibitem{GiZai2}
A. Giambruno and M. Zaicev. 
\emph{Central polynomials of associative algebras and their growth}. 
Proc. Amer. Math. Soc. 147 (2019) 909--919.

\bibitem{GiambrunoZaicev} A. Giambruno and M. Zaicev.  {\em Exponential codimension growth of PI algebras: an exact estimate.} Adv. Math. 142 (1999) 221--243.

\bibitem{GZ14}
A. Giambruno and M. Zaicev.
\emph{ Growth of polynomial identities: is the sequence of codimensions eventually non-decreasing?}.
Bull. Lond. Math. Soc. {46} (2014), 771--778.


\bibitem{GiambrunoZai} A. Giambruno and M. Zaicev. Polynomial Identities and Asymptotic Methods. Mathematical Surveys and Monographs, vol. 122.
American Mathematical Society, Providence, RI, 2005.

  \bibitem{Sandra} S. M. A. Jorge and A. C. Vieira. {\em On minimal varieties of quadratic growth.} Linear Algebra Appl. 418 (2006) 925--938.

  \bibitem{Kemer}
A. R. Kemer.
\emph{Ideals of Identities of Associative Algebras}.
Translations of Mathematical Monographs, vol.~87.
American Mathematical Society, Providence, RI, 1991.

      \bibitem{Kem} A. R. Kemer. \textit{ Varieties of finite rank}. In: Proc. 15th All-Union Algebraic Conf., vol. 2, Krasnoyarsk, 1979, p. 73 (in Russian).

  \bibitem{PKS} P. Koshlukov, A. Krasilnikov and
E. A. da Silva. \textit{The central polynomials for the finite dimensional Grassmann algebras.}
Algebra Discrete Math. no. 3 (2009)  69--76.

\bibitem{KoshlukovLaMattina2015}
P. Koshlukov and D. La Mattina.
{\it Graded algebras with polynomial growth of their codimensions.}
J. Algebra {434} (2015) 115--137.

  \bibitem{Daniela} D. La Mattina.  
		\newblock{\em Characterizing varieties of colength $\leq 4$.} \newblock Comm. Algebra 37 (2009) 1793--1807.

  \bibitem{Danielasub}  D. La Mattina. \newblock{\em Varieties of almost polynomial growth: classifying their subvarieties}. \newblock Manuscripta Math. 123 (2007) 185--203.

    \bibitem{RG} A. Regev.  {\em Existence of identities in $A\otimes B$}.   Israel J. Math. {11} (1972) 131--152.

	
\end{thebibliography}
\end{document}